\documentclass[fontsize=10pt]{scrreprt}

\usepackage[square,numbers,nonamebreak]{natbib}
\usepackage{setspace}

\usepackage[english]{babel}

\usepackage{chngcntr}

\usepackage[hang]{footmisc}
 at 11truept
\font\ssc=pplrc9d at 11 truept
\newcommand\qedbox{$\rlap{$\sqcap$}\sqcup$}

\usepackage[hmarginratio=1:1, bottom=2cm, top=2.5cm]{geometry}

\usepackage[automark]{scrlayer-scrpage}
\cohead{}

\let\ceheadL\cehead
\renewcommand\cehead[1]{
\ceheadL{\textnormal{#1}}
}

\cehead{F. de Giovanni -- I. de las Heras -- M. Trombetti}
\lefoot{}
\rofoot{}

\usepackage{graphicx}
\usepackage[usenames,dvipsnames,svgnames,table]{xcolor}
\definecolor{Maroon}{cmyk}{0, 0.87, 0.68, 0.32}
\definecolor{RoyalBlue2}{cmyk}{80,100,0,0.1}

\newcommand\auths[1]{\large \textsc{\textcolor{Maroon}{#1}}\setstretch{1.2}}
\newcommand\titl[1]{\center \linespread{1.1}\color{RoyalBlue2}\Large\textbf{ #1}\color{black}\bigskip} 
\renewcommand\abstract[1]{
\begin{center}
{\textbf{Abstract}}
\end{center}
{
\linespread{1.1}\fontsize{9pt}{-10pt}\selectfont #1}}

\usepackage[utf8]{inputenc}
\usepackage[sc,osf]{mathpazo}
\usepackage[euler-digits]{eulervm}
\usepackage[T1]{fontenc}
\usepackage{mathtools}

\DeclareSymbolFont{operators}{\encodingdefault}{ppl}{m}{n}
\DeclareMathAlphabet{\mathbf}{\encodingdefault}{ppl}{bx}{n}
\DeclareMathAlphabet{\mathit}{\encodingdefault}{ppl}{m}{it}

\usepackage{amsfonts}
\usepackage{amsmath}
\usepackage{ragged2e}
\usepackage{titlesec}
\usepackage{amssymb}

\renewcommand{\thesection}{\arabic{section}}
 
\titleformat{\section}{\medskip\bigskip\normalfont\Large\bf}{\thesection}{0.5em}{}
\titleformat{\subsection}{\smallskip\bigskip\normalfont\large\bf}{\thesubsection}{0.5em}{}

\usepackage{amsthm}
\newtheoremstyle{dotless}{}{}{\itshape}{}{\bfseries}{}{1em}{}

\newtheorem*{theo*}{Theorem}
\newtheorem*{rem*}{Remark}
\theoremstyle{dotless}
\newtheorem{theo}{Theorem}

\newtheorem{lem}[theo]{Lemma}
\newtheorem{cor}[theo]{Corollary}

\newtheorem{quest}[theo]{Question}
\renewenvironment{proof}{\smallbreak\noindent {\sc Proof \;---\;}}{\hfill\qedbox}

\counterwithout{theo}{section}
\numberwithin{theo}{section}

\DeclareOldFontCommand{\rm}{\normalfont\rmfamily}{\mathrm}
\DeclareOldFontCommand{\sf}{\normalfont\sffamily}{\mathsf}
\DeclareOldFontCommand{\tt}{\normalfont\ttfamily}{\mathtt}
\DeclareOldFontCommand{\bf}{\normalfont\bfseries}{\mathbf}
\DeclareOldFontCommand{\it}{\normalfont\itshape}{\mathit}
\DeclareOldFontCommand{\sl}{\normalfont\slshape}{\@nomath\sl}
\DeclareOldFontCommand{\sc}{\normalfont\scshape}{\@nomath\sc}

\usepackage{soul}
\DeclareSymbolFont{newfont}{OML}{cmm}{m}{it}%
\DeclareMathSymbol{\Varrho}{3}{newfont}{37}

\usepackage{tikz}

\begin{document}

\titl{On the Lattice of Closed Subgroups of a Profinite Group\footnote{The first and third authors are supported by GNSAGA (INdAM) and are members of the non-profit association ``Advances in Group Theory and Applications'' (www.advgrouptheory.com).
The second author is supported by the Spanish Government project PID2020-117281GB-I00,  
partially by FEDER funds and by the Basque Government project IT483-22.
He has also received funding from the European Union’s Horizon 2021 research and innovation programme under the Marie Sklodowska-Curie grant agreement, project 101067088.\\
On behalf of all authors, the corresponding author states that there is no conflict of interest.\\
This manuscript has no associated data.}}

\auths{Francesco de Giovanni -- Iker de las Heras -- Marco Trombetti}

\thispagestyle{empty}
\justify\noindent
\setstretch{0.3}
\abstract{The subgroup lattice of a group is a great source of information about the structure of the group itself. The aim of this paper is to use a similar tool for studying profinite groups. In more detail, we study the lattices of closed or open subgroups of a profinite group and its relation with the whole group. We show, for example, that procyclic groups are the only profinite groups with a distributive lattice of closed or open subgroups, and we give a sharp characterization of profinite groups whose lattice of closed (or open) subgroups satisfies the~Dedekind modular law; we actually give a precise description of the behaviour of modular elements of the lattice of closed subgroups. We also deal with the problem of carrying some structural information from a profinite group to another one having an isomorphic lattice of closed (or open) subgroups. Some interesting consequences and related results concerning decomposability and the number of profinite groups with a given lattice of closed  (or open) subgroups are also obtained.}

\setstretch{2.1}
\noindent
{\fontsize{10pt}{-10pt}\selectfont {\it Mathematics Subject Classification \textnormal(2020\textnormal)}: 20E18, 20E15}\\[-0.8cm]

\noindent 
\fontsize{10pt}{-10pt}\selectfont  {\it Keywords}: profinite group; closed subgroup lattice; modular lattice; distributive lattice; decomposable lattice; width of a lattice\\[-0.8cm]

\section{Introduction}
\setstretch{1.1}
\fontsize{11pt}{12pt}\selectfont

{\it Profinite groups} are topological groups that are compact, Hausdorff and totally disconnected. The interest of this kind of  groups arose first from Galois theory, and in fact it turns out that profinite groups are precisely the Galois groups endowed with the~Krull topology (see Section 2.11 of \cite{Zalesski} for more details).
In the last decades, the abstract and topological structure of profinite groups has been extensively studied by many authors (see for instan\-ce~\cite{Acciarri}, \cite{topft}, \cite{KhuShu}, \cite{Klopsch}, \cite{Wilson}); we refer to the monograph \cite{Zalesski} for results and terminology concerning profinite groups.

The aim of this paper is to introduce a new point of view in studying profinite groups: we will give descriptions of a profinite group starting from the way the closed (or open) subgroups are ``positioned''. For instance, two of our main results (see Theorems \ref{distributive} and \ref{modular}) describe what a profinite group looks like when it does not contain closed subgroups $A,B,C,D,E$ positioned as follows 
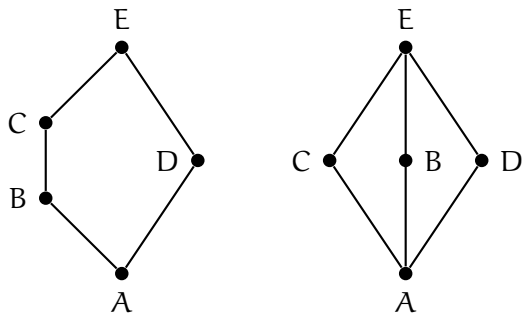
\begin{center}
\begin{tikzpicture}[thick, circ/.style={shape=circle, fill=black, inner sep=1pt, scale=1.5, draw, node contents=}]
\draw node (c1) at (1, 3) [circ, label=above:{$E$}];
\draw node (c2) at (0, 2) [circ, label=left:{$C$}];
\draw node (c3) at (0, 1) [circ, label=left:{$B$}];
\draw node (c4) at (2, 1.5) [circ, label=left:{$D$}];
\draw node (c5) at (1, 0) [circ, label=below:{$A$}];
\draw (c1) to (c2) to (c3) to (c5) to (c4) to (c1);
\end{tikzpicture} $\quad\quad$
\begin{tikzpicture}[thick, circ/.style={shape=circle, fill=black, inner sep=1pt, scale=1.5, draw, node contents=}]
\draw node (c1) at (1, 3) [circ, label=above:{$E$}];
\draw node (c2) at (0, 1.5) [circ, label=left:{$C$}];
\draw node (c3) at (2, 1.5) [circ, label=right:{$D$}];
\draw node (c4) at (1, 0) [circ, label=below:{$A$}];
\draw node (c5) at (1, 1.5) [circ, label=right:{$B$}];
\draw (c1) to (c2) to (c4) to (c3) to (c1) to (c4);
\end{tikzpicture}
\end{center}
In the above pictures, the lines describe "adjacent" lattice relations between closed subgroups, meaning that there is no proper closed subgroup between the subgroups in the endpoints of the lines; thus, for example, in both pictures, the smallest closed subgroup containing~$B$ and~$D$ is $E$, while $A$ is the intersection of $C$ and $D$.

\smallskip

Applications of lattice theory to abstract algebra dates back to 1928, when E. Fischer noted that the fundamental theorem of Galois theory simply stated the existence of an isomorphism between the subgroup lattice of a Galois group and a certain lattice of fields. He consequently proposed the problem of determining whether two groups with the same subgroup lattice are isomorphic or not. One of his students, Ada Rottländer, answered this question in the negative \cite{ada}, but this paved the way to the lattice-theoretical isomorphism problem in algebra, an instance of which is: classify all groups that share the same subgroup lattice. Thus, for example, locally cyclic groups can be characterized as those groups whose subgroup lattice is distributive (see \cite{ore}), and distributivity is easily detectable in a lattice, as it reduces to the absence of pentagonal and diamond sublattices. Another example is given by the Iwasawa characterization of groups in which the Dedekind modular law holds (see \cite{Schmidt}, Chapter~2). The study of the subgroup lattice has many related interesting consequences, and in fact it gives useful conditions under which a group is decomposable as a non-trivial free product (see \cite{Schmidt}, Theorem 7.1.20). Our hope is that the study of the closed (or open) subgroup lattice of a profinite group could led to similar conditions characterizing non-trivial free products of profinite groups. Having this in mind, we wish this paper to be a reference point for future work on the subject, and in fact, we also highlight many interesting problems in this area, which, if solved, would provide a great deal of information about profinite groups themselves.

\smallskip





We start in Section \ref{Sec: Lattices} by giving a detailed account of some general properties of the lattice of closed (or open) subgroups of a profinite group.
In particular we show that if a topologically finitely generated profinite group $G$ has the same subgroup lattices as a profinite group $H$, then they must also have the same close (or open) subgroup lattice, and $H$ must be topologically finitely generated too (see Theorem \ref{topfg}).

Section 3 is devoted to the study of profinite groups whose closed or open subgroup lattice is distributive. We will show that these profinite groups are precisely the procyclic ones (see~The\-o\-rem~\ref{distributive}), and that this result has some interesting consequences concerning profinite groups having isomorphic open subgroup lattices (see for instance The\-o\-rem~\ref{thopenfg}).


It is much more involved dealing with profinite groups whose lattice of closed subgroups has some restrictions on the cardinality of maximal antichains, or with modular elements of the lattice of closed subgroups of a profinite group. 
We will do so in Sections~\ref{width} and~\ref{mod}, respectively.

Profinite groups with a decomposable lattice of closed subgroups are dealt with in~Section~\ref{decomposability}.


\medskip

Concerning terminology, since we deal with both group theoretical concepts and their topological analogous, we adhere to the rule that whenever we work with a topological concept, we add the adjective ``topological''. For example, if a subgroup~$H$ of a profinite group $G$ is a ``topological Sylow subgroup of $G$'', then this must be understood in the sense explained in~\cite{Zalesski}; if we only say that~$H$ is a ``Sylow subgroup of~$G$'', we mean that $H$ is a maximal element of the set of all $p$-subgroups of $G$, for some prime $p$. 
We refer to \cite{Rob72} as a general reference for concepts and notation concerning abstract group theory, but we point out that: 

\begin{itemize}
    \item[(1)] $\mathbb{Z}_p$ denotes the additive group of $p$-adic integers (endowed with its profinite topology).
    
    \item[(2)] $\mathcal{C}_n$ denotes the cyclic group of order $n$.
    
    \item[(3)] If $G$ is a profinite group, then $N\trianglelefteq_oG$ means that $N$ is an open normal subgroup of $G$.
    
    \item[(4)] If $Y\leq X$ are subgroups of a group $G$, then $[X/Y]$ denotes the lattice made by all subgroups of $G$ lying between $X$ and $Y$. Moreover, if $G$ is profinite and $Y\leq X$ are closed subgroups of $G$, we put $[X/Y]_c$ to be the set of all closed subgroups of~$G$ lying between $X$ and $Y$. Note also that, if $Y\leq X$ are open subgroups of the profinite group $G$, then the interval $[X/Y]$ is made by open subgroups and so there is no need to introduce the symbol $[X/Y]_o$.
    
    \item[(5)] If $\mathcal{X}$ is any family of (profinite) groups, the unrestricted direct product of the members of this family will be denoted by $\operatorname{Cr}_{X\in\mathcal{X}}X$
    , while the restricted direct product will be denoted by~\hbox{$\operatorname{Dr}_{X\in\mathcal{X}}X$}.
    
    \item[(6)] If $G$ is a group (resp. profinite group), we denote by $\pi(G)$ (resp. $\pi^\star(G)$) the set of all primes $p$ for which $G$ has non-trivial $p$-elements (resp. has non-trivial procyclic pro-$p$ subgroups); two periodic abstract groups (resp. profinite groups)~$X$ and $Y$ are said to be {\it coprime} if $\pi(X)\cap\pi(Y)=\emptyset$ (resp. \hbox{$\pi^\star(X)\cap\pi^\star(Y)=\emptyset$)}; of course, if $X$ and $Y$ are coprime subgroups of a (profinite) group, then $X\cap Y=\{1\}$.
    
    \item[(7)] Topological closures will always be denoted by a bar and referred to simply as closures, so if $H$ is an abstract subgroup of the profinite group $G$, its closure in $G$ will be denoted by $\overline{H}$.
    
    
    \item[(8)] If $X$ is any subgroup of a group $G$, then the normal core of $X$ in $G$ will be denoted by $X_G$, while the normal closure of $X$ in $G$ by $X^G$.
\end{itemize}


\section{Lattices in a profinite group}
\label{Sec: Lattices}

Let $G$ be a profinite group. The set of all closed subgroups of $G$ will be denoted by~$\mathcal{L}_c(G)$, so if $X$ is any subgroup of $G$, the closure $\overline{X}$ of $X$ in $G$ belongs to~$\mathcal{L}_c(G)$. If $X,Y\in\mathcal{L}_c(G)$, we define the {\it closed join} $X\,{\mbox{\footnotesize$\vee$}}_cY$ of $X$ and $Y$ as $\overline{X\,{\mbox{\footnotesize$\vee$}}\, Y}$, where, as usual,~\hbox{$X\,{\mbox{\footnotesize$\vee$}}\,Y=\langle X,Y\rangle$.} Thus, using the symbol ${\mbox{\footnotesize$\wedge$}}$ to denote the intersection of elements in $\mathcal{L}_c(G)$, we have that
$\big(\mathcal{L}_c(G),{\mbox{\footnotesize$\vee$}_c}\,,{\mbox{\footnotesize$\wedge$}}\big)$ is a lattice that will be referred to as the {\it closed subgroup lattice} of $G$; in particular, if $Y\leq X$ are closed subgroups of $G$, the set~$[X/Y]_c$ is naturally a sublattice of $\mathcal{L}_c(G)$, and if $N$ is a closed normal subgroup of~$G$, then~$\mathcal{L}_c(G/N)\simeq[G/N]_c$.

Note that~$\mathcal{L}_c(G)$ is not in general a sublattice of the lattice~$\mathcal{L}(G)$ of all subgroups of $G$. To see this, it is enough to consider the semidirect product~\hbox{$\langle x\rangle\ltimes\mathbb{Z}_p$,} where $p$ is any prime and $x$ is the inverting automorphism of~$\mathbb{Z}_p$; here, the subgroup generated by any conjugate of~$x$ has order~$2$, so is closed, but the subgroup generated by any two distinct subgroups of order~$2$ is countably infinite, so it cannot be closed. The situation is somewhat better in the universe of abelian profinite groups. In fact, by~Ty\-chonoff's theorem, the join of two permuting closed subgroups of any profinite group is still closed, and so in particular $\mathcal{L}_c(G)$ is a sublattice of $\mathcal{L}(G)$ when~$G$ is abelian; in this context, note that the consideration of an unrestricted direct product of infinitely many $\mathbb{Z}_p$ (for a given prime~$p$) shows that~$\mathcal{L}_c(G)$ need not be a complete sublattice of~$\mathcal{L}(G)$, even for abelian profinite groups. It will follow from our~The\-o\-rem~\ref{modular} that unrestricted direct products of coprime profinite groups with a modular closed subgroup lattice form another class of profinite groups~$G$ for which the closed subgroup lattice is a sublattice of the lattice of all subgroups.

The situation is much better if we consider the set $\mathcal{L}_o(G)$ of all open subgroups of~$G$: it is in fact easy to see that this is always a sublattice of both $\mathcal{L}(G)$ and $\mathcal{L}_c(G)$, although not always a complete sublattice (if $G$ is infinite, $\mathcal{L}_o(G)$ does not have a minimum). This lattice will be referred to as the {\it open subgroup lattice} of $G$. It turns out that $\mathcal{L}_o(G)$ is completely determined by $\mathcal{L}_c(G)$, and in fact the following easy lemma shows that open subgroups of profinite groups can be recognized within the closed subgroup lattice.

\begin{lem}\label{openclosed}
Let $G$ be a profinite group and let $X$ be a closed subgroup of $G$. Then $X$ is open if and only if it is contained in only  finitely many closed subgroups of $G$.
\end{lem}
\begin{proof}
Let $X$ be an open subgroup of $G$. Then $X$ is a closed subgroup of finite index of $G$, and hence the set ${\cal L}_X$ of all closed subgroups of $G$ containing~$X$ is finite. Suppose conversely that ${\cal L}_X$ is finite, so that in particular there are only finitely many open subgroups containing $X$.
Since $X$ is the intersection of the open subgroups of $G$ containing it, it follows that $X$ is open in $G$.
\end{proof}

\medskip

Since the closed subgroups of finite index are precisely the open subgroups of a profinite group, the above result should be seen in relation with a theorem of~Za\-cher~\cite{Zacher} and Rips stating that the finiteness of the index of a subgroup in an abstract group can be recognized from the lattice of its subgroups. An obvious consequence of Lemma \ref{openclosed} is that the finiteness of a profinite group can be recognized by any of the lattices $\mathcal{L}(G)$, $\mathcal{L}_c(G)$, and $\mathcal{L}_o(G)$.

\begin{cor}\label{corfinito}
Let $G$ be a profinite group. The following conditions are equivalent:
\begin{itemize}
    \item[\textnormal{(1)}] $G$ is finite;
    \item[\textnormal{(2)}] $\mathcal{L}(G)$ is finite;
    \item[\textnormal{(3)}] $\mathcal{L}_c(G)$ is finite;
    \item[\textnormal{(4)}] $\mathcal{L}_o(G)$ is finite.
\end{itemize}
\end{cor}

\begin{quest}\label{constructlcg}
Let $G$ be a profinite group. Can we construct $\mathcal{L}_c(G)$ starting from $\mathcal{L}_o(G)$?
\end{quest}

Let $G$ and $H$ be profinite groups. Recall that a {\it projectivity} from $G$ onto $H$ is a lattice isomorphism from $\mathcal{L}(G)$ onto $\mathcal{L}(H)$. Similarly, we define a {\it $c$-projectivity} as a lattice isomorphism from $\mathcal{L}_c(G)$ onto $\mathcal{L}_c(H)$, and an {\it $o$-projectivity} as a lattice isomorphism from $\mathcal{L}_o(G)$ onto $\mathcal{L}_o(H)$. Of course, it follows from Lemma \ref{openclosed} that every~\hbox{$c$-pro}\-jectivity induces an~\hbox{$o$-pro}\-jectivity. Question \ref{constructlcg} can be essentially rephrased as follows in terms of projectivities. 

\begin{quest}\label{constructlcg2}
Is it true that any $o$-projectivity can be extended to a $c$-projectivity?
\end{quest}

Although our aim is to show that the lattices $\mathcal{L}_c(G)$ and $\mathcal{L}_o(G)$ have a strong impact on the structure of a profinite group $G$, it should be noted that none of the latti\-ces~\hbox{$\mathcal{L}_c(G)$, $\mathcal{L}_o(G)$, $\mathcal{L}(G)$} is able to characterize profinite groups in the universe of topological groups. This can be seen from the following examples:

\begin{itemize}
    \item[\textnormal{(1)}] If $p$ is any prime, then $\mathbb{Z}_p$ endowed with the usual profinite topology obviously has the same subgroup lattice of $\mathbb{Z}_p$ with the discrete topology. The same remark applies to any infinite profinite group.
    \item[\textnormal{(2)}] Let $p$ be a prime and let $G=\mathbb{Z}_p$  with the usual profinite topology. Consider a topological generator $x$ of $G_1$ and let $H=\langle x\rangle$ endowed with the subspace topology induced by $G$. Since $H$ is Hausdorff and countable, it is not compact. However, $\mathcal{L}_c(G)\simeq \mathcal{L}_c(H)$ and $\mathcal{L}_o(G)\simeq \mathcal{L}_o(H)$.
\end{itemize}

\begin{quest}
Let $G$ and $H$ be profinite groups with isomorphic subgroup lattice. Is it true that $\mathcal{L}_c(G)\simeq\mathcal{L}_c(H)$? Or, is it at least true that $\mathcal{L}_o(G)\simeq\mathcal{L}_o(H)$?
\end{quest}

The above question can be answered at least in the finitely generated case. Recall that a profinite group is {\it topologically finitely generated} if it is the closure of a finitely generated subgroup. In the case of topologically finitely generated profinite groups, a well known theorem of Novikov and Segal \cite{Niksegal} states that every subgroup of finite index is open, so the open subgroups are precisely the subgroups of finite index.

\begin{theo}\label{topfg}
Let $G$ and $H$ be profinite groups such that $\mathcal{L}(G)\simeq\mathcal{L}(H)$. If $G$ is topologically finitely generated, then so is $H$. 
Moreover, we have $\mathcal{L}_c(G)\simeq\mathcal{L}_c(H)$ and $\mathcal{L}_o(G)\simeq\mathcal{L}_o(H)$.
\end{theo}
\begin{proof}
Let $\varphi$ be a projectivity from $G$ onto $H$, and let $x_1,\ldots,x_n$ be the topological generators of $G$.
Write $X_i=\langle x_i\rangle$ for $i=1,\ldots,n$, and observe that $X_i^{\varphi}$ is a cyclic subgroup of $H$.
Define
$$
V=\overline{\langle X_1^{\varphi},\ldots,X_n^{\varphi}\rangle}.
$$
Since $V$ contains the subgroups $X_i^{\varphi}$, the subgroup $U=V^{\varphi^{-1}}$ contains the subgroups~$X_i$.
Moreover, since $V$ is closed in $H$, it is the intersection of the open subgroups containing it, so in particular it is the intersection of finite index subgroups of~$H$.
Therefore, $U$ is intersection of finite index subgroups of $G$, and since $G$ is finitely generated, if follows that $U$ is closed in $G$.
In particular, $U=G$, so \hbox{$V=U^{\varphi}=G^{\varphi}=H$} is topologically finitely generated.

Now, since the finiteness of the index can be recognized in the subgroup lattice, and the open subgroups of $G$ and of $H$ are precisely those of finite index, it follows that~$\varphi$ induces a lattice isomorphism from $\mathcal{L}_o(G)$ onto $\mathcal{L}_o(H)$.
Moreover, as every closed subgroup is the intersection of the open subgroups containing it, it also follows that~$\varphi$ induces a lattice isomorphism from $\mathcal{L}_c(G)$ to $\mathcal{L}_c(H)$, and the theorem follows.
\end{proof}

\medskip

The previous statement should be compared to Corollary \ref{thclosedfg} and Theorem \ref{thopenfg}. Moreover, the proof of Theorem \ref{topfg} shows that any projectivity $\varphi$ of $G$ onto $H$ induces a~\hbox{$c$-pro}\-jectivity; the following lemma shows that in such cases $\varphi$ preserves closures of subgroups.

\begin{lem}
Let $\varphi$ be a projectivity from the profinite group $G$ onto the profinite group $H$. If~$\varphi$ induces a $c$-projectivity between $G$ and $H$, then $\overline{X}^{\varphi}=\overline{X^{\varphi}}$ for every $X\leq G$.
\end{lem}
\begin{proof}
Since $X\le\overline{X}$, we have $X^{\varphi}\le\overline{X}^{\varphi}$, and since $\overline{X}^{\varphi}$ is closed in $H$, we obtain $\overline{X^{\varphi}}\le \overline{X}^{\varphi}$.
From $X=(X^{\varphi})^{\varphi^{-1}}$, a symmetric argument shows that $\overline{X}\le\overline{X^{\varphi}}^{\varphi^{-1}}$, so in particular $\overline{X}^{\varphi}\le\overline{X^{\varphi}}$, and the statement follows.
\end{proof}

\medskip

In the forthcoming sections we essentially describe profinite groups whose closed subgroup lattice (resp. open subgroup lattice) satisfies one of the following conditions: distributivity (Section \ref{distributivity}), decomposability (Section \ref{decomposability}), existence of a finite bound on the width (Section \ref{width}) and modularity (Section \ref{mod}).

\section{Distributivity}\label{distributivity}

In this section, we deal with profinite groups whose closed or open subgroup lattice is distributive. Recall that a lattice $(L,{\mbox{\footnotesize$\vee$}},{\mbox{\footnotesize$\wedge$}})$ is {\it distributive} if one of the following equivalent identities 
\begin{itemize}
    \item[(1)] $x\,{\mbox{\footnotesize$\vee$}}\,(y\,{\mbox{\footnotesize$\wedge$}}\,z)=(x\,{\mbox{\footnotesize$\vee$}}\,y)\,{\mbox{\footnotesize$\wedge$}}\,(x\,{\mbox{\footnotesize$\vee$}}\,z)$ for all $x,y,z\in L$
    \item[(2)] $x\,{\mbox{\footnotesize$\wedge$}}\,(y\,{\mbox{\footnotesize$\vee$}}\,z)=(x\,{\mbox{\footnotesize$\wedge$}}\,y)\,{\mbox{\footnotesize$\vee$}}\,(x\,{\mbox{\footnotesize$\wedge$}}\,z)$ for all $x,y,z\in L$
\end{itemize} holds. The main result for distributivity (see Theorem \ref{distributive}) is not too difficult to prove and its proof exploits the fact mentioned in the introduction that the subgroup lattice of a group is distributive if and only if every finitely generated subgroup is cyclic. Nevertheless, Theorem \ref{distributive} has some interesting corollaries concerning the finite topological generation (see Corollary \ref{thclosedfg} and Theorem \ref{thopenfg}).

\begin{theo}\label{distributive}
Let $G$ be a profinite group. The following statements are equivalent:
\begin{itemize}
    \item[\textnormal{(1)}] $\mathcal{L}_c(G)$ is distributive;
    \item[\textnormal{(2)}] $\mathcal{L}_o(G)$ is distributive;
    \item[\textnormal{(3)}] $G$ is procyclic.
\end{itemize}
\end{theo}
\begin{proof}
Of course, (1) implies (2). Suppose now $\mathcal{L}_o(G)$ is distributive. If $N$ is any open normal subgroup of $G$, then~$[G/N]$ is distributive, so~$G/N$ is cyclic by Co\-rol\-la\-ry~1.2.4 of \cite{Schmidt}.
Thus, $G$ is procyclic, so (2) implies (3).

Assume finally that $G$ is procyclic. Then $G$ is an unrestricted direct product of groups of coprime order of type $\mathbb{Z}_p$ or ${\cal C}_{q^n}$, where $p$ and $q$ are primes.
Now, the closed subgroup lattice of these groups are chains, and chains are obviously distributive.
Moreover, the closed subgroup lattice of the unrestricted direct product of coprime groups is the unrestricted direct product of the closed subgroup lattices of the groups.
Thus, since distributivity is preserved by taking unrestricted direct products, it follows that $\mathcal{L}_c(G)$ is distributive, and hence (3) implies (1).
\end{proof}

\medskip
Theorem \ref{distributive} in particular says that, if $G$ is a profinite group and $X$ is a closed subgroup of $G$, we can understand whether $X$ is procyclic or not just by looking at $\mathcal{L}_c(G)$. This means that the following analog of Lemma \ref{topfg} holds.

\begin{cor}\label{thclosedfg}
Let $G$ be a topologically finitely generated profinite group.
Then the image of $G$ under any $c$-projectivity is topologically finitely generated as well.
\end{cor}

The corresponding analog for $o$-projectivities is much more subtle. For its proof we need to use the properties of modular elements of a lattice (see \cite{Schmidt}, Theorem 2.1.5); recall that an element $m$ of the lattice $(L,{\mbox{\footnotesize$\vee$}},{\mbox{\footnotesize$\wedge$}})$ is {\it modular} in $L$ if it satisfies the following two conditions: \begin{itemize}
    \item[(1)] $x\,{\mbox{\footnotesize$\vee$}}\,(m\,{\mbox{\footnotesize$\wedge$}}\, z)=(x\,{\mbox{\footnotesize$\vee$}}\, m)\,{\mbox{\footnotesize$\wedge$}}\, z$ for all $x,z\in L$ with $x\leq z$,
    \item[(2)] $m\,{\mbox{\footnotesize$\vee$}}\,(y\,{\mbox{\footnotesize$\wedge$}}\, z)=(m\,{\mbox{\footnotesize$\vee$}}\, y)\,{\mbox{\footnotesize$\wedge$}}\, z$ for all $y,z\in L$ with $m\leq z$.
\end{itemize}

It is clear that the image of any open (resp. closed) normal subgroup of a profinite group under any $o$-projectivity (resp. $c$-projectivity) is still a modular element of the open subgroup lattice (resp. closed subgroup lattice).

\begin{lem}
    \label{lem: filteredbelow}
    Let $G$ be a profinite group and $H$ a closed subgroup of $G$.
    Let $\{U_i\mid i\in I\}$ be a family of closed subgroups of $G$ filtered from below, and write $U=\bigcap_{i\in I} U_i$.
    If $UH$ is open in $G$, then there exists $i\in I$ such that $UH=U_iH$.
\end{lem}
\begin{proof}
    By Proposition 2.1.4 (a) of~\cite{Zalesski}, we have $UH=\bigcap_{i\in I}(U_iH).$
    Now, since $UH$ has finite index in $G$ and since $\{U_i\mid i\in I\}$ is filtered from below, there exists $i$ such that $UH=U_iH$, as desired.
\end{proof}

\begin{theo}\label{thopenfg}
The image of any topologically finitely generated profinite group under an~\hbox{o-pro}\-jectivity is topologically finitely generated as well.
\end{theo}
\begin{proof}
Let $\varphi$ be an $o$-projectivity from a topologically finitely generated profinite group $G$ onto a profinite group $H$. Let $x_1,\ldots,x_n$ be topological generators of $G$, and write $X_i=\overline{\langle x_i\rangle}$ for all $i$. Then $$X_i=\bigcap_{Y\in\mathcal{O}_i}Y$$ where $\mathcal{O}_i$ is the set of all open subgroups of $G$ containing $X_i$. Put now $$\mathcal{O}_i^\varphi=\bigl\{Y^\varphi\;|\; Y\in\mathcal{O}_i\bigr\}$$ and $$Z_i=\bigcap_{Z\in\mathcal{O}_i^\varphi}Z.$$
Of course, $Z_i^\varphi$ is a closed (although not necessarily open) subgroup of $H$ and the set~$\mathcal{O}_i^\varphi$ contains the intersection of any pair of its members. 

Let us show that $Z_i$ is procyclic. To this aim, it is enough to prove that $Z_i/(Z_i\cap N)$ is cyclic for every $N\trianglelefteq_o H$.
Fix such an $N$ and let $M\trianglelefteq_oG$ be such that $M\leq N^{\varphi^{-1}}$. By Lemma \ref{lem: filteredbelow}, we can find $Y\in\mathcal{O}_i$ such that
$$
Y^\varphi N=Z_i N\quad\textnormal{and}\quad YM=X_iM.
$$
Therefore, the subgroup lattice of $Z_i/(Z_i\cap N)\simeq Y^\varphi N/N$ is isomorphic to the interval $$\big[(Y\,{\mbox{\footnotesize$\vee$}}\,N^{\varphi^{-1}})/\,N^{\varphi^{-1}}\big].$$ On the other hand, $N^{\varphi^{-1}}$ is a modular element of $\mathcal{L}_o(G)$, so $$\big[(Y\,{\mbox{\footnotesize$\vee$}}\,N^{\varphi^{-1}})/\,N^{\varphi^{-1}}\big]\simeq\big[Y/(Y\cap N^{\varphi^{-1}})\big]$$ by Theorem 2.1.5 of \cite{Schmidt}. But the latter interval is a sublattice of the distributive lattice $$[Y/(Y\cap M)]\simeq[YM/M]=\mathcal{L}\bigl(X_i/(X_i\cap M)\bigr),$$ so it is distributive as well.
Thus, the subgroup lattice of $Z_i/(Z_i\cap N)$ is distributive and $Z_i$ is procyclic.

It remains to show that the closed join of all $Z_i$ is the whole group $H$. Again, it suffices to see that $\langle X_i^\varphi\,|\, i=1,\ldots,n\rangle N/N=H/N$ for all $N\trianglelefteq_oH$. As before, fix $N\trianglelefteq_o H$ and let~\hbox{$M\trianglelefteq_oG$} be such that $M\leq N^{\varphi^{-1}}$.
By Lemma \ref{lem: filteredbelow}, we may find open subgroups $Y_1,\ldots,Y_n$ of $G$ such that
$$
Y_i^\varphi N=Z_i N\quad\textnormal{and}\quad Y_iM=X_iM
$$
for all $i=1,\ldots,n$. Now,
$$
\langle Z_i\,|\, i=1,\ldots,n\rangle N/N=\langle Y_i^\varphi\,|\, i=1,\ldots,n\rangle N/N.
$$
On the other hand,
$$\Big(\langle Y_i^\varphi\,|\, i=1,\ldots,n\rangle N\Big)^{\varphi^{-1}}
=\bigvee_{i=1,\ldots,n}Y_i\,{\mbox{\footnotesize$\vee$}}\,N^{\varphi^{-1}}\geq \bigvee_{i=1,\ldots,n}Y_iM=\bigvee_{i=1,\ldots,n}X_iM=G$$ and consequently
$$
\langle Z_i\,|\, i=1,\ldots,n\rangle N/N=H/N.
$$
Thus $\overline{\langle Z_i\,|\, i=1,\ldots,n\rangle}=H$ and we are done.
\end{proof}

\medskip

Actually, the arguments of the above proof can be used to say something more: let $\varphi$ be an $o$-projectivity from the profinite group $G$ onto the profinite group~$H$, let $X$ be a closed subgroup of $G$, and let $\mathcal{O}$ be the set of all open subgroups of $G$ containing $X$, so that $X$ is the intersection of all members of $\mathcal{O}$. If we choose~$X^\varphi$ as the intersection of all $Y^\varphi$ with $Y\in\mathcal{O}$, we have defined a map from~$\mathcal{L}_c(G)$ to~$\mathcal{L}_c(H)$. The proof (and the arguments) of Theorem \ref{thopenfg} shows that this map has the following properties:
\begin{itemize}
    \item[1)] If $X$ is a procyclic subgroup of $G$, then $X^\varphi$ is procyclic.
    \item[2)] If $X$ is a closed subgroup of $G$ which is topologically generated by a family $\mathcal{F}$ of closed subgroups, then $X^\varphi$ is topologically generated by the images of the subgroups in $\mathcal{F}$ under $\varphi$ (in other words, $\varphi$ is a complete homomorphism with respect to $\mbox{\footnotesize$\vee$}_c$). In particular, if $Y\leq X$ are closed subgroups of $G$, then $Y^\varphi\leq X^\varphi$.
\end{itemize}

We call this map the {\it $c$-extension} of $\varphi$, and we usually denote it with $\varphi$ too. Thus in the following sections, if $\varphi : \mathcal{L}_o(G)\longrightarrow\mathcal{L}_o(H)$ is an~\hbox{$o$-pro}\-jectivity and $X$ is a closed subgroup of $G$, we denote by $X^\varphi$ the image of $X$ by the $c$-extension of $\varphi$.
Thus, we can further rephrase Questions \ref{constructlcg} and \ref{constructlcg2} in the following way.

\begin{quest}\label{projectext}
Let $\varphi$ be an $o$-projectivity between two profinite groups~$G$ and $H$. Is the~\hbox{$c$-exten}\-sion of $\varphi$ a $c$-projectivity between $G$ and $H$?
\end{quest}

Note that a positive answer to the above question is equivalent to proving that $c$-extensions are bijective (see for instance~\cite{Schmidt},~The\-o\-rem~1.1.2). Some weaker (but reasonable) forms of Question \ref{projectext} are the following ones.

\begin{quest}
Let $\varphi$ be an $o$-projectivity between two profinite groups~$G$ and $H$, and let~$X$ be a closed subgroup of $G$. Is it true that $\mathcal{L}_o(X)\simeq \mathcal{L}_o(X^\varphi)$? 
\end{quest}

\begin{quest}
Let $\varphi$ be an $o$-projectivity between two profinite groups~$G$ and $H$. Is the~\hbox{$c$-exten}\-sion of $\varphi$ a homomorphism with respect to $\mbox{\footnotesize$\wedge$}$\,?
\end{quest}

\section{Decomposability}\label{decomposability}

A lattice $L$ is said to be {\it directly decomposable} if it is isomorphic to the direct product of two non-trivial lattices. Groups whose subgroup lattice is directly decomposable have a clear structure: they are precisely the direct products of two non-trivial coprime periodic groups. More generally, there is a correspondence between unrestricted direct decompositions of the subgroup lattice of a group $G$ and restricted direct decompositions of $G$ in terms of coprime periodic groups (see for instance \cite{Schmidt}, The\-o\-rem~1.6.5).

The main result of this section completely solves the corresponding problem for profinite groups: it turns out that the description is very similar once one replaces restricted direct products by unrestricted ones. However, as we shall see after the proof, this is one of the few cases in which there is a little bit of asymmetry between the behaviour of the closed subgroup lattice and that of the open subgroup lattice.
Finally, we will show in the end of the section that direct decomposability can be detected by only looking at the quotient modulo the topological Frattini subgroup (see~The\-o\-rem~\ref{thtopfrat}).

\begin{theo}\label{suzuki}
Let $G$ be a profinite group and let $\{L_i\}_{i\in I}$ be a family of lattices with $|L_i|\geq2$ for all $i\in I$. The following statements are equivalent:
\begin{itemize}
    \item[\textnormal{(1)}] $\mathcal{L}_c(G)\simeq\operatorname{Cr}_{i\in I}L_i$\,;
    \item[\textnormal{(2)}] $G\simeq\operatorname{Cr}_{i\in I}G_i$, where $G_i$ is a profinite group with $\mathcal{L}_c(G_i)\simeq L_i$ for every $i\in I$, and $\pi^\star(G_i)\cap\pi^\star(G_j)=\emptyset$ for every $i,j\in I$ with~\hbox{$i\neq j$}.
\end{itemize}
\end{theo}
\begin{proof}
Clearly, (2) implies (1) since any closed subgroup $X$ of $G$ decomposes as~\hbox{$\operatorname{Cr}_{i\in I}X_i$,} where $X_i\leq G_i$, and consequently $$\mathcal{L}_c(G)\simeq\underset{i\in I}{\operatorname{Cr}}\,\mathcal{L}_c(G_i).$$

Conversely, assume (1) and write $L=\operatorname{Cr}_{i\in I}L_i$. Let $\sigma:\mathcal{L}_c(G)\longrightarrow L$ be a lattice isomorphism.
For each $i\in I$, let $0_j$ and $1_j$ be, respectively, the minimum and maximum of $L_i$, and define
$$
F_i=\underset{j\in I}{\operatorname{Cr}}\,X_j\le L,
$$
where $X_i=1_i$ and $X_j=0_j$ for $j\neq i$.
Put \hbox{$G_i=\sigma^{-1}(F_i)$}.
Of course, $G_i\in\mathcal{L}_c(G)$, $\mathcal{L}_c(G_i)\simeq L_i\simeq[f_i/0]$, and $$\overline{\langle G_j\,:\, j\in I\rangle}=G.$$

If $N\trianglelefteq_o G$, then
$$
\mathcal{L}_c(G/N)\simeq[1/N^\sigma]\simeq
\underset{i\in I}{\operatorname{Cr}}\,\big[F_i\,{\mbox{\footnotesize$\vee$}}\,N^\sigma/N^\sigma\big].
$$
Since $F_i\,{\mbox{\footnotesize$\vee$}}\,N^\sigma$ corresponds to $G_iN$ in the previous sequence of isomorphisms, we may use Theorem 1.6.5 of \cite{Schmidt} to get $\pi(G_iN/N)\cap\pi(G_jN/N)=\emptyset$ and $[G_i,G_j]\leq N$ for all $i\neq j$. The arbitrariness of $N$ yields
$$
\pi^\star(G_i)\cap\pi^\star(G_j)=\emptyset\quad\textnormal{and}\quad [G_i,G_j]=\{1\}
$$
for all $i\neq j$, and hence $G\simeq\operatorname{Cr}_{i\in I}G_i$. The proof is complete.
\end{proof}

\begin{cor}
Let $G$ be a profinite group. Then $\mathcal{L}_c(G)$ is directly decomposable if and only if~\hbox{$G=G_1\!\times\! G_2$,} where $G_1$ and $G_2$ are non-trivial profinite groups with \hbox{$\pi^\star(G_1)\cap\pi^\star(G_2)=\emptyset$.}
\end{cor}

\begin{cor}
Let $G$ be a profinite group. Then $\mathcal{L}(G)$ is directly decomposable if and only if~$G$ is periodic and $\mathcal{L}_c(G)$ is directly decomposable.
\end{cor}
\begin{proof}
If $\mathcal{L}(G)$ is directly decomposable, then $G=G_1\times G_2$, where $G_1$ and $G_2$ are periodic subgroups with $\pi(G_1)\cap\pi(G_2)=\emptyset$ (see  \cite{Schmidt}, The\-o\-rem~1.6.5). In particular, $G_1$ and $G_2$ are closed and~\hbox{$\mathcal{L}_c(G)\simeq\mathcal{L}_c(G_1)\times\mathcal{L}_c(G_2)$} by Theorem \ref{suzuki}. The conserve follows similarly.
\end{proof}

\medskip

We note that Theorem \ref{suzuki} does not hold in its generality if we replace $\mathcal{L}_c(G)$ by~\hbox{$\mathcal{L}_o(G)$.} In fact, an infinite unrestricted direct product of non-trivial lattices cannot be the open subgroup lattice of any profinite group. Thus, for example, the subgroup $\operatorname{Cr}_{p\in\mathbb{P}}p\mathbb{Z}_p$ is not an open subgroup of $\operatorname{Cr}_{p\in\mathbb{P}}(\mathbb{Z}_p)$. On the other hand, if we only deal with finite direct products, then we get the following analog of Theorem \ref{suzuki}. Note that, as usual, we write $0$ and $1$ for the eventual minimum and maximum elements of a lattice.

\begin{theo}\label{suzuki2}
Let $G$ be a profinite group and let $\{L_i\}_{i\in I}$ be a finite family of lattices with \hbox{$|L_i|\geq2$} for all $i\in I$. The following statements are equivalent:
\begin{itemize}
    \item[\textnormal{(1)}] $\mathcal{L}_o(G)\simeq\operatorname{Cr}_{i\in I}L_i$;
    \item[\textnormal{(2)}] $G\simeq\operatorname{Cr}_{i\in I}G_i$, where $G_i$ is a profinite group with $\mathcal{L}_o(G_i)\simeq L_i$ for every $i\in I$, and $\pi^\star(G_i)\cap\pi^\star(G_j)=\emptyset$ for every $i,j\in I$ with~\hbox{$i\neq j$}.
\end{itemize}
\end{theo}
\begin{proof}
Suppose first that (2) holds. Since $I$ is finite, every open subgroup $X$ of~$G$ decomposes as $\operatorname{Cr}_{i\in I}X_i$, where $X_i\leq G_i$, and hence (1) follows.

Conversely, in order to prove that (1) implies (2), we may assume $\mathcal{L}_o(G)\simeq K\times L$ for only two non-trivial lattices $K$ and $L$. Let $\sigma:\mathcal{L}_o(G)\longrightarrow K\times L$ be an isomorphism, and put $$G_K=\bigcap_{\ell\in L}(1,\ell)^{\sigma^{-1}}\quad\textnormal{and}\quad G_L=\bigcap_{k\in K} (k,1)^{\sigma^{-1}}.$$ In particular, $G_K$ and $G_L$ are closed subgroups of $G$.
If $N\trianglelefteq_oG$, by Lemma \ref{lem: filteredbelow} we have
$$G_KN=(1,\ell_K)^{\sigma^{-1}}N\quad\textnormal{and}\quad G_LN=(k_L,1)^{\sigma^{-1}}N$$ for suitable $\ell_K\in L$ and $k_L\in K$ with $N^\sigma\geq (k_L,\ell_K)$. Now, the subgroup lattice of $G/N$ is the direct product of the lattices
$$
\big[(1,\ell_K)\,\mbox{\footnotesize$\vee$}\,N^\sigma/N^\sigma\big]\quad\textnormal{and}\quad\big[(k_L,1)\,\mbox{\footnotesize$\vee$}\,N^\sigma/N^\sigma\big].
$$
This means that $$G/N=G_KN/N\times G_LN/N$$ and $\pi(G_KN/N)\cap\pi(G_LN/N)=\emptyset$ by Theorem 1.6.9 of \cite{Schmidt}. The arbitrariness of $N$ yields that $G=G_K\times G_L$ and that $\pi^\star(G_K)\cap\pi^\star(G_L)=\emptyset$.

We need finally to show that $\mathcal{L}_o(G_K)\simeq K$ and $\mathcal{L}_o(G_L)\simeq L$.
Let $U$ be any open subgroup of $G_K$.
Again, since $UG_L=U(a,1)$ for some $a\in K$ by Lemma \ref{lem: filteredbelow}, there exists $a_U\in K$ such that
$$
(U\times G_L\big)^\sigma=(a_U,1).
$$
Conversely, let $a_U\in K$ and suppose $(a_U,1)^{\sigma^{-1}}=U\times V$, for some open subgroups~$U$ of $G_K$ and $V$ of $G_L$.
Then, if $N\trianglelefteq_oG$, we have $VN/N=G_LN/N$, which means that $V=G_L$.

This shows that the map $$\psi: U\in\mathcal{L}_o(G_K)\longmapsto a_U\in K$$ is a lattice isomorphism. Similarly, $\mathcal{L}_o(G_L)\simeq L$, and the proof is complete.~\end{proof}

\begin{cor}
Let $G$ be a profinite group. Then $\mathcal{L}_c(G)$ is directly decomposable if and only if $\mathcal{L}_o(G)$ is directly decomposable.
\end{cor}

As a consequence of the proofs of Theorems \ref{suzuki} and \ref{suzuki2}, we have the following.

\begin{cor}
Let $G=\operatorname{Cr}_{i\in I}G_i$ be a profinite group, where $\pi^\star(G_i)\cap\pi^\star(G_j)=\emptyset$ for every $i,j\in I$ with~\hbox{$i\neq j$.} If $\varphi$ is a $c$-projectivity from $G$ onto $H$, then $H=\operatorname{Cr}_{i\in I} G_i^\varphi$ and \hbox{$\pi^\star(G_i^\varphi)\cap\pi^\star(G_j^\varphi)=\emptyset$} if $i\neq j$. The same conclusion holds if $I$ is finite and $\varphi$ is an~\hbox{$o$-pro}\-jectivity.
\end{cor}

\begin{cor}
Let $\pi$ be a set of primes and let $X$ and $Y$ be closed subgroups of a profinite group~$G$ such that $X$ is generated by pro-$\pi$ subgroups, $Y$ by pro-$\pi'$ subgroups and~\hbox{$[X,Y]=\{1\}$.} If $\varphi:G\longrightarrow H$ is an $o$-projectivity, then~\hbox{$[X^\varphi,Y^\varphi]=\{1\}$.} Moreover, if~\hbox{$\pi^\star(X)\cap\pi^\star(Y)=\emptyset$,} then $\pi^\star(X^\varphi)\cap\pi^\star(Y^\varphi)=\emptyset$.
\end{cor}
\begin{proof}
Suppose first $\pi^\star(X)\cap\pi^\star(Y)=\emptyset$.
The argument we employed in the first half of the proof of Theorem \ref{thopenfg} shows that $\pi^\star(X^\varphi)\cap\pi^\star(Y^\varphi)=\emptyset$ and \hbox{$[X^\varphi,Y^\varphi]=\{1\}$.}


Now, if $U$ is a pro-$\pi$ subgroup of $X$ and $V$ is a pro-$\pi'$ subgroup of $Y$, it follows that~\hbox{$[U^\varphi,V^\varphi]=\{1\}$.} Since $X^\varphi$ (resp. $Y^\varphi$) is generated by the images of the pro-$\pi$ subgroups of $X$ (resp. of $Y$)
, we obtain that $[X^\varphi,Y^\varphi]=\{1\}$. 
\end{proof}

\medskip

Note that in the previous statement, the properties $[X,Y]=\{1\}$ and \hbox{$\pi^\star(X)\cap\pi^\star(Y)=\emptyset$} need not be preserved individually under projectivities. To see this, it is enough to consider the projectivity between the symmetric group of degree $3$ and the elementary abelian group of order~$9$.

\smallskip

In the final part of this section, we show that one can understand if the closed subgroup lattice (or the open subgroup lattice) of a profinite group is directly decomposable by looking at the quotient modulo the topological Frattini subgroup. Recall that the {\it topological Frattini subgroup} $\Phi_c(G)$ of a profinite group $G$ is the intersection of all maximal open subgroups of $G$ (see \cite{Zalesski}, Section 2.8).

\begin{lem}\label{quozientefrattini}
Let $G$ be a profinite group. If $p\in\pi^\star(G)$, then $p\in\pi^\star\big(G/\Phi_c(G)\big)$.
\end{lem}
\begin{proof}
Let $p\in\pi^\star(G)$ and choose $N\trianglelefteq_oG$ such that $p\in\pi(G/N)$.  Then~$p$ divides the order of $(G/N)/\Phi(G/N)$, where $\Phi(G/N)$ is the Frattini subgroup of the finite group $G/N$ (see for instance \cite{Rob82}, p.263). Now, Corollary 2.8.3 of \cite{Zalesski} yields that $p$ belongs to~$\pi^\star\big(G/\Phi_c(G)\big)$ and proves the statement.
\end{proof}

\begin{theo}\label{thtopfrat}
Let $G$ be a profinite group. The following properties are equivalent:
\begin{itemize}
    \item[\textnormal{(1)}] $\mathcal{L}_c(G)$ is directly decomposable;
    \item[\textnormal{(2)}] $\mathcal{L}_o(G)$ is directly decomposable;
    \item[\textnormal{(3)}] There exist non-trivial closed subgroups $H$ and $K$ of $G$ such that $G=H\times K$ and\linebreak $\pi^\star(H)\cap\pi^\star(K)=\emptyset$;
    \item[\textnormal{(4)}] $\mathcal{L}_c\big(G/\Phi_c(G)\big)$ is directly decomposable;
    \item[\textnormal{(5)}] $\mathcal{L}_o\big(G/\Phi_c(G)\big)$ is directly decomposable.
\end{itemize}
\end{theo}
\begin{proof}
Statements (1), (2) and (3) are all equivalent by Theorem \ref{suzuki}. Similarly,~(4) and (5) are equivalent by the same result. Thus, it suffices to prove that~(3) and (4) are equivalent.

Suppose (3) holds. Then $$G/\Phi_c(G)=H\Phi_c(G)/\Phi_c(G)\times K\Phi_c(G)/\Phi_c(G),$$
where the factors $H\Phi_c(G)/\Phi_c(G)$ and~\hbox{$K\Phi_c(G)/\Phi_c(G)$} are non-trivial and coprime.

Conversely, suppose (4) holds and write $$G/\Phi_c(G)=H_1/\Phi_c(G)\times H_2/\Phi_c(G),$$ where $\pi^\star\big(H_1/\Phi_c(G)\big)\cap\pi^\star\big(H_2/\Phi_c(G)\big)=\emptyset$. Corollary 2.8.4 of \cite{Zalesski} yields that $\Phi_c(G)$ is pronilpotent, so it is the direct product of its topological Sylow subgroups.
Let~$D_i$ be the subgroup generated by those topological Sylow $p$-subgroups of $\Phi_c(G)$ such that $p\in\pi^\star\big(H_i/\Phi_c(G)\big)$; in particular, $D_i\trianglelefteq G$, and $\Phi_c(G)=D_1\times D_2$ since every prime~$p$ in~$\pi^\star\big(\Phi_c(G)\big)$ is also contained in~$\pi^\star\big(G/\Phi_c(G)\big)$ by Lemma \ref{quozientefrattini}.
Now, $D_1$ is a topological Hall subgroup of $H_2$, so by the Schur--Zassenhaus theorem (see \cite{Zalesski},~The\-o\-rem~2.3.15), there is a (non-trivial) closed subgroup $G_2$ of $H_2$ such that $H_2=G_2\ltimes D_1$.
By the Frattini argument, we have $$G=N_G(G_2)H_2=N_G(G_2)D_1,$$
and since $D_1\leq\Phi_c(G)$, we obtain $G=N_G(G_2)$.
Thus, $G_2$ is normal in $G$.
Similarly, we define $G_1$ and show that it is normal in $G$. Moreover, $G_i\geq D_i$ and $$G=H_1H_2=G_1D_2G_2D_1=G_1G_2.$$ Since clearly $\pi^\star(G_1)\cap\pi^\star(G_2)=\emptyset$, we have that $G=G_1\times G_2$ and (3) is proved.
\end{proof}

\section{Finite width}\label{width}

Let $G$ be a finite group whose subgroup lattice has a non-trivial chain as a direct factor. Then there are infinitely many non-isomorphic finite groups having the same subgroup lattice of~$G$. Conversely, it turns out that if $L$ is any finite lattice with no chain as a direct factor, then there are only finitely many groups whose subgroup lattice is isomorphic to $L$ (see \cite{Schmidt}, Theorem 1.6.10).

Moving from finite to infinite groups, these results are no longer true, even when we deal with the closed (or open) subgroup lattice of a profinite group.
To see this, we first construct a counterexample for abstract infinite groups and adapt it later to profinite groups.
Let $G_1$ be a group of order $2$ and, for each integer~\hbox{$i>1$,} let~$G_i$ be an elementary abelian group of order $p_i^2$, where~$p_i$ is the $i$th prime number.
Let~$G$ be the restricted direct product of all $G_i$'s. Clearly, the subgroup lattice of $G$ has a non-trivial chain as a direct factor, but if~$H$ is any group such that~\hbox{$\mathcal{L}(H)\simeq\mathcal{L}({G})$,} then~\hbox{$H\simeq G$.}
Indeed, $H$ must be the restricted direct product of coprime finite subgroups $H_i$ which are lattice-isomorphic to the subgroups~$G_i$ (see Theorem 2.2.5 of \cite{Schmidt}), so that $H_1$ has prime order and, for $i>1$, either $H_i\simeq G_i$ or $H_i\simeq \mathcal{C}_{q_i}\ltimes\mathcal{C}_{p_i}$ for some prime~\hbox{$q_i\neq p_i$.} Thus~$p_i$ divides $|H_i|$ for each $i>1$. It follows that $|H_1|=2$ and so also $H_i\simeq G_i$ for all~$i$, because the subgroups $H_i$ have coprime orders. Therefore $G$ is isomorphic to $H$.

Now, in a similar way, if $G^\ast$ is the unrestricted direct product of all $G_i$'s, we see, using Theorem \ref{suzuki}, that if $H$ is any profinite group with $\mathcal{L}_o(H)\simeq\mathcal{L}_o(G^\ast)$, then~\hbox{$H\simeq G^\ast$.}

This shows that there exist infinite groups (resp. infinite profinite groups) whose subgroup lattice (resp. closed subgroup lattice, open subgroup lattice) has a chain as a direct factor, but which are uniquely determined by such a lattice, so the first result on finite groups we mentioned at the beginning of the section cannot be extended to the infinite case. Now, we prove that also the second-mentioned result does not hold for infinite (profinite) groups. 

In the case of abstract groups, it would be enough to consider {\it Tarski $p$-groups} (that is, infinite groups whose proper subgroups are cyclic of prime order $p$), but for the sake of generality we now describe a construction that, as done before, can easily be adapted to the profinite case.
By Dirichlet's theorem on arithmetic progressions, we can find an infinite set $\mathcal P$ of triples of pairwise distinct prime numbers such that:

\begin{itemize}
    \item[1)] if $(p_1,p_2,p_3)\neq(q_1,q_2,q_3)$ are elements of $\mathcal P$, then $\{p_1,p_2,p_3\}\cap\{q_1,q_2,q_3\}=\emptyset$;
    \item[2)] if $(p_1,p_2,p_3)\in\mathcal P$, then $p_2p_3$ divides $p_1-1$.
\end{itemize}

\noindent
For each $\mathbf p=(p_1,p_2,p_3)\in\mathcal P$ and $i=2,3$, let $G_{\mathbf p,i}=\mathcal{C}_{p_i}\ltimes \mathcal{C}_{p_1}$ be the only \hbox{non-abelian} group of order $p_1p_i$. Moreover, let $\mathcal F$ be the set of all functions $f$ from $\mathcal P$ into $\{2,3\}$, and for each~\hbox{$f\in\mathcal F$,} let $G_f$ be the restricted direct product of the groups $G_{\mathbf p,f(\mathbf p)}$ with~\hbox{$\mathbf p\in\mathcal P$.} It is easy to see that, for all $f,g\in\mathcal{F}$, the lattice $\mathcal{L}(G_f)$ has no chains as direct factors, and~\hbox{$\mathcal{L}(G_f)\simeq\mathcal{L}(G_g)$.} Thus, there are $2^{\aleph_0}$-many non-isomorphic groups whose subgroup lattices are isomorphic and have no chains as direct factors (see Theorem 2.2.5 of \cite{Schmidt}).

Similarly, if $G_f^\ast$ is the unrestricted direct product of the groups $G_{\mathbf p,f(\mathbf p)}$ endowed with the natural profinite topology, then $\mathcal{L}_c(G_f^\ast)\simeq\mathcal{L}_c(G_g^\ast)$ and $\mathcal{L}_o(G_f^\ast)\simeq\mathcal{L}_o(G_g^\ast)$ for all~\hbox{$f,g\in\mathcal F$.} Again, we have found $2^{\aleph_0}$-many non-isomorphic profinite groups whose closed subgroup lattices (resp. open subgroup lattices) are isomorphic and have no chains as direct factors; we use here the results proved in Section \ref{distributivity}.

\smallskip
Nevertheless, something can still be achieved in this direction. In the abstract setting, if the lattice~$L$ has finite {\it width} $w$ (that is, the antichains of $L$ are finite and $w$ is the maximum order of an antichain), then groups whose subgroup lattice is isomorphic to $L$ can be classified (see \cite{Brandl}, where the width is called {\it Dilworth number}). Here we prove that a classification is possible also in the profinite case. First, we need to recall that a~\hbox{\it $BFC$-group} is a group having a uniform finite bound on the conjugacy classes of its elements; a celebrated theorem of B.H. Neumann states that a group is $BFC$ if and only if its derived subgroup is finite.

\begin{theo}\label{thwidth}
Let $G$ be an infinite profinite group. The following conditions are equivalent:
\begin{itemize}
    \item[\textnormal{(1)}] $\mathcal{L}_c(G)$ has finite width;
    \item[\textnormal{(2)}] $\mathcal{L}_o(G)$ has finite width;
    \item[\textnormal{(3)}] $G=K\ltimes T$, where $K\simeq\mathbb{Z}_p$ and $T$ is a finite $p'$-group  for some prime number $p$.
\end{itemize}
\end{theo}
\begin{proof}
The fact that (1) implies (2) is obvious. Conversely, assume (2) and let $w$ be the width of $\mathcal{L}_o(G)$.
If $X_1,\ldots,X_{n}$ are closed subgroups of $G$ such that  $X_i\not\leq X_j$ for all~\hbox{$i,j\in\{1,\ldots,n$\}} with $i\neq j$, then we can find an open normal subgroup $N$ of $G$ such that $X_iN\not\leq X_jN$ for all $i,j\in\{1,\ldots,n\}$ with $i\neq j$. Since $X_iN$ is open for every $i$, we get that $n\leq w$, so the width of $\mathcal{L}_c(G)$ is $w$ as well. Therefore (1) and (2) are equivalent.

Thus, it suffices to show that (1) and (3) are equivalent. Assume (1) and let $w$ be the width of~$\mathcal{L}_c(G)$.
Clearly, $|\pi^*(G)|\le w$ and, by Corollary 1 of \cite{Brandl}, the topological~Sy\-low subgroups of $G$ are either infinite procyclic or finite.

We first deal with the case in which $G$ is abelian.
Let $p$ be a prime such that the~Sy\-low~\hbox{$p$-sub}\-group $K$ of $G$ satisfies $K\simeq \mathbb{Z}_p$, and note that for any prime $q\neq p$, the group $\mathbb{Z}_p\times \mathbb{Z}_q$ has infinite width.
This shows that $K$ is the unique infinite Sylow subgroup of $G$, and therefore $G=K\times T$, where $T$ is the Hall $p'$-subgroup of $G$, as desired.

We now deal with the general case.
We first claim that $G$ is a $BFC$-group.
For the primes $q\in\pi^*(G)$ for which the Sylow $q$-subgroups of $G$ are finite, we write $n_q$ for the order of such subgroups.
We will show that $|G:C_G(x)|$ is uniformly bounded in terms of $w$, $\pi^*(G)$ and the $n_q$'s for every $x\in G$.

Fix $x\in G$ and observe that the conjugacy classes of~$X=\overline{\langle x\rangle}$ form an antichain of~$\mathcal{L}_c(G)$, so $|G:N_G(X)|\le w$.
If $X$ is finite, then $|X|$ is bounded in terms of the $n_r$'s.
Thus, $|\text{Aut}(X)|$ is also bounded in terms of the $n_r$'s, and hence so is $|N_G(X):C_G(X)|$.
We may therefore assume that $X$ is infinite.
Now, for every $q\in \pi^*(G)$, let $X_q$ be the Sylow $q$-subgroup of $X$, and let $p$ be the unique prime in $\pi^*(G)$ for which $X_p$ is infinite (recall that $X$ is procyclic, and hence abelian).
Let $G_p$ be a Sylow $p$-subgroup of $G$ containing $X_p$, and note that $G_p$ is also procyclic by Corollary 1 of \cite{Brandl}.
For $q\neq p$, observe that $|G_p:N_{G_p}(X_q)|\leq w$ and that $|N_{G_p}(X_q):C_{G_p}(X_q)|$ is bounded in terms of $n_q$.
Hence, for $C=\bigcap_{r\neq p} C_{G_p}(X_r)$, we have that $|G_p:C|$ is bounded in terms of~$w$ and of the $n_r$'s.
Moreover, $C\le C_G(X)$, so in particular the Sylow $p$-subgroup of~$N_G(X)/C_G(X)$ is bounded in terms of $w$ and the $n_r$'s. It is easy to see that the
rest of the (topological) Sylow subgroups of $N_G(X)/C_G(X)$ are bounded in terms of $w$, $p$, and the $n_q$'s, so, since $N_G(X)/C_G(X)$ is abelian, the claim follows.

Therefore $G$ is a $BFC$-group, and consequently $G'$ is finite.
Now,
$$
G/G'\simeq P/G'\times T/G',
$$
where $T/G'$ is finite and $P/G'\simeq\mathbb{Z}_p$, for some prime $p$. Thus, $p$ is the unique prime in $\pi^*(G)$ for which the Sylow $p$-subgroups of $G$ are infinite (and hence procyclic).
Therefore,~\hbox{$G=K\ltimes T$,} where~$K\simeq\mathbb{Z}_p$
and $T$ is a finite $p'$-group, so (3) is proved.

\smallskip

Finally, we assume that $G$ satisfies (3), and that $\mathcal{L}_c(G)$ has infinite width. Thus, for any $n\in\mathbb{N}$, we may find an antichain $\mathcal{A}_n$ of $\mathcal{L}_c(G)$ having cardinality $n$. Since $T$ is finite, we can certainly assume that every element of $\mathcal{A}_n$ is not contained in $T$, and consequently that every element of $\mathcal{A}_n$ is infinite.

Let $U$ be any infinite closed subgroup of $G$. Then $U=Y_U\ltimes(U\cap T)$, where $Y_U\simeq\mathbb{Z}_p$, and $|G:U|$ is finite. Since $T$ is finite, the antichains $\mathcal{A}_n$ can be chosen in such a way that for any $U\in\mathcal{A}_n$ and $V\in\mathcal{A}_m$, we have $U\cap T=V\cap T$. Moreover, any $Y_U$, $U\in\mathcal{A}_n$, is certainly contained in a conjugate of $K$ (being $K$ a topological Sylow $p$-subgroup of~$G$). But there are finitely many such conjugates, and hence the antichains $\mathcal{A}_n$ can be chosen in such a way that there exists $g\in G$ such that for any $U\in\mathcal{A}_n$ and $V\in\mathcal{A}_m$, the closed subgroups $Y_U$ and $Y_V$ are contained in $K^g$. This is clearly a contradiction, because then the elements of $\mathcal{A}_n$ are comparable. This proves (1) and completes the proof of the theorem.
\end{proof}

\begin{cor}
Let $L$ be a lattice with finite width $w$. Then there are at most countably many non-isomorphic profinite groups whose closed \textnormal(resp. open\textnormal) subgroup lattice is isomorphic to $L$.
If, moreover, $L$ has no non-trivial chain as a direct factor, then there are at most finitely many non-isomorphic profinite groups whose closed \textnormal(resp. open\textnormal) subgroup lattice is isomorphic to $L$.
\end{cor}
\begin{proof}
By Theorem 1.6.10 of \cite{Schmidt} and Corollary \ref{corfinito}, we may assume $L$ is infinite. Consider now a profinite group $G$ whose closed (resp. open) subgroup lattice is isomorphic to $L$. Then $G$ has the structure described in (3) of Theorem \ref{thwidth}, and of course there are at most countably many non-isomorphic groups of this type.

Suppose now $L$ has no non-trivial chain as a direct factor, and
write $G=K\ltimes T$, whe\-re~\hbox{$K\simeq \mathbb{Z}_p$} and $T$ a finite $p'$-group for some prime number $p$.
Let $C=C_K(T)$, so~\hbox{$C<K$,} and define $H=(K/C)\ltimes T$.

Observe that $H$ is finite and $\mathcal{L}(H)$ is (isomorphic to) a fixed finite sublattice of~$L$. In fact, it is easy to see that $C$ is the largest infinite procyclic  pro-$p$ subgroup~of the centre of $G$, and so its image in $L$ can be recognized. Moreover, $\mathcal{L}(H)$ cannot have any non-trivial chain as a direct factor, as otherwise also~$L$ would. In particular, by~Theorem 1.6.10 of \cite{Schmidt}, there are finitely many possibilities for $p$ and $T$, and the proof is complete.
\end{proof}

\medskip

Note that if we have some particular bounds on the width of the~open subgroup lattice of a profinite group, additional information on its structure can be obtained by using for example Proposition 5 of \cite{Brandl}.

\section{Modularity}\label{mod}

In this section we characterize profinite groups whose closed subgroup lattice satisfies the Dedekind modular law and we deal in more detail with modular elements of the closed subgroup lattice. Recall first that a lattice $(L,{\mbox{\footnotesize$\vee$}},{\mbox{\footnotesize$\wedge$}})$ is said to be {\it modular} (or to satisfy the {\it Dedekind modular law}) if
$$
x\,{\mbox{\footnotesize$\vee$}}\,(y\,{\mbox{\footnotesize$\wedge$}}\,z)=(x\,{\mbox{\footnotesize$\vee$}}\,y)\,{\mbox{\footnotesize$\wedge$}}\,z
$$
for all $x,y,z\in L$ with $x\leq z$.
Modularity is easily detectable in a lattice since one just needs to check that there are no pentagonal sublattices (see \cite{Schmidt}, Theorem 2.1.2). Obviously, a lattice is modular if and only if all its elements are modular.
Relevant examples of groups whose subgroup lattice is modular are those in which every pair of subgroups~$H$ and $K$ is {\it permutable} (that is, $HK=KH$); in particular, all abelian groups and more generally all Dedekind groups have a modular subgroup lattice. The structure of {\it modular} groups (that is, groups with a modular subgroup lattice) has been completely determined by Iwasawa in a series of papers (see \cite{Schmidt} for further details). In particular, for our purposes, we should keep in mind the following results: 

\begin{itemize}
    \item[(1)] Theorem 2.3.1 of \cite{Schmidt}, which gives a detailed description of finite $p$-groups with a modular subgroup lattice. This result makes use of the so-called {\it Iwasawa triples} $(A,b,s)$ for a $p$-group $G$ (see \cite{Schmidt}, p.60): here, $A$ is an abelian subgroup of $G$, $b\in G$, $s$ is a positive integer which is at least $2$ in case $p=2$,  $G=A\langle b\rangle$, and \hbox{$b^{-1}ab=a^{1+p^s}$} for all $a\in A$. More in general, (non-abelian) locally finite $p$-groups with a modular subgroup lattice are described in Theorem 2.4.14 of \cite{Schmidt}. It turns out that these groups are precisely those $p$-groups in which any two subgroups permute (see \cite{Schmidt}, Lemma 2.3.2).
    \item[(2)] The\-o\-rem~2.4.13 of \cite{Schmidt}, which states that a locally finite group is modular if and only if it is a direct product of $P^*$-groups and modular $p$-groups with coprime orders. Recall that a~{\it $P^*$-group} is just a semidirect product of an elementary abelian normal subgroup $A$ by a cyclic group $\langle t\rangle$ of prime power order such that~$t$ induces a {\it power automorphism} of prime order on $A$ (see \cite{Schmidt}, p.69), that is, an automorphism of prime order fixing all subgroups of $A$. 
    \item[(3)] It easily turns out from (1) and (2) that every locally finite group which is locally modular must be modular.

\end{itemize}

Before characterizing profinite groups with a modular closed or open subgroup lattice, we need the following easy lemma.

\begin{lem}\label{permute}
Let $G$ be a profinite group, and let $H$ and $K$ be closed subgroups of $G$. If $HKN=KHN$ for every open normal subgroup $N$ of $G$, then $HK=KH$.
\end{lem}
\begin{proof}
By Tychonoff's theorem, the products $HK$ and $KH$ are closed subsets of~$G$. Therefore~$HK$ (resp. $KH$) is the intersection of all subgroups of the form $HKN$ (resp.~$KHN$), where~$N$ is open and normal in $G$. By hypothesis $HKN=KHN$ for all such $N$, and hence $HK=KH$.
\end{proof}

\begin{theo}\label{modular}
Let $G$ be a profinite group. The following statements are equivalent:
\begin{itemize}
    \item[\textnormal{(1)}] $\mathcal{L}_c(G)$ is modular;
    \item[\textnormal{(2)}] $\mathcal{L}_o(G)$ is modular;
    \item[\textnormal{(3)}] $G$ is an unrestricted direct product of coprime profinite groups each of which is of one of the following types: \begin{itemize}
        \item[\textnormal{(i)}] abelian;
        \item[\textnormal{(ii)}] locally finite and modular;
        \item[\textnormal{(iii)}] $X\ltimes A$, where $X\simeq\mathbb{Z}_p$ for some prime $p$, $A$ is an elementary abelian $q$-group for some prime $q\neq p$ and $X=\overline{\langle x\rangle}$, where $x$ acts on $A$ as a power automorphism of prime order;
        \item[\textnormal{(iv)}] $X\ltimes A$, where $X\simeq\mathbb{Z}_p$ for some prime $p$, $A$ is an abelian pro-$p$ subgroup, and there are $k\in\mathbb{N}$ and $x\in X$  such that $X=\overline{\langle x\rangle}$ and $a^x=a^{1+p^k}$ for every $a\in A$; moreover, $k\geq2$ if $p=2$.
    \end{itemize}
\end{itemize}
\noindent If any of the above equivalent statements holds, then the abstract subgroup generated by two arbitrary closed subgroups of $G$ is closed, and, if $G$ is of type \textnormal{(iv)}, every pair of closed subgroups of $G$ is permutable.
\end{theo}
\begin{proof}
Clearly, we may assume $G$ is infinite and non-abelian.

It is immediate that (1) implies (2). Assume (2). If $N$ is any open normal subgroup of $G$, then $G/N$ is modular, so it is the direct product of~\hbox{$P^*$-groups} and modular $p$-groups of coprime orders (see \cite{Schmidt}, The\-o\-rem~2.4.4).

Let $\pi=\pi^*(G)$. Let $\sigma$ be the partition of~$\pi$ satisfying the following property: two distinct primes $p,q\in\pi$ are in the same element of $\sigma$ if and only if there is an open normal subgroup $N_{p,q}$ of $G$ such that~$G/N_{p,q}$ contains a~\hbox{$P^*$-group} whose order is divisible by $pq$. Note that every element of $\sigma$ has cardinality at most $2$.

Let now $M$ and $N$ be open normal subgroups of $G$ such that $M\le N$.
For every~\hbox{$\tau\in \sigma$,} the $\tau$-components of~$G/M$ are mapped onto the $\tau$-components of $G/N$ under the natural epimorphism $G/M\longrightarrow G/N$.
This means that for each~\hbox{$\tau\in\sigma$,} the topological $\tau$-component of $G$ is unique and $G$ splits into the unrestricted direct product of these components.

In order to prove (3), it is therefore possible to assume that $G$ coincides with one of these components. Assume first $G$ is the inverse limit of finite $P^*$-groups whose order is a $\{p,q\}$-number, where $p<q$. Then $G$ has an elementary abelian normal~\hbox{$q$-sub}\-group~$A$ such that $A$ is closed and $G/A\simeq\mathbb{Z}_p$.
By the Schur--Zassenhaus theorem, we have $G=X\ltimes A$, for some procyclic pro-$p$ subgroup $X$. Let $x$ be a topological generator of $X$. Then, in any homomorphic image of $G$ by an open normal subgroup $N$, $xN$ acts on $AN/N$ as a power automorphism of prime order. This means that $x$ acts on $A$ as a power automorphism of prime order. If $X$ is finite we are thus in case (ii), otherwise we are in case (iii).

Suppose now $G$ is the inverse limit of finite modular $p$-groups, so $G$ is a pro-$p$ group. Thus, if $N$ is any open normal subgroup of $G$, then $G/N$ satisfies~The\-o\-rem~2.3.1 of~\cite{Schmidt}. If every such quotient is Dedekind, then $G$ is a Dedekind $p$-group and we are done. We may hence assume every such quotient satisfies point (b) of The\-o\-rem~2.3.1. In particular, we see that $\overline{G'}$ is abelian. 

Assume first $G/\overline{G'}$ is periodic. If $\overline{G'}$ is periodic, then $G$ is periodic, so also locally finite; it follows that $G$ is locally a modular $p$-group, which means that $G$ is modular and locally finite.
It is therefore safe to assume that $\overline{G'}$ is not periodic, so then $\overline{G'}$ contains an infinite procyclic subgroup~$V$.
Now,~The\-o\-rem~2.3.1~(b) of \cite{Schmidt} entails that~$V$ is normal in $G$ and, moreover, that $V/V^4$ is central in $G$ if $p=2$.
Since $V$ has no automorphism whose order is a power of $p$ acting like this, and since~\hbox{$V\leq \overline{G'}\leq C_G(V)$,} we have that $$C_G(V)=N_G(V)=G.$$
Thus, $V$ is central in $G$. But $\overline{G'}$ is generated by its infinite procyclic subgroups, and hence $\overline{G'}$ is contained in the centre of $G$. Now, $G/\overline{G'}$ has finite exponent~(see~\cite{Zalesski},~Lem\-ma 4.3.7) and consequently a theorem of Mann (see \cite{mann}) yields that $G'$ has finite exponent, so~$\overline{G'}$ has finite exponent, a contradiction.

Suppose $G/\overline{G'}$ is non-periodic. We assign to every open normal subgroup $M$ of $G$ an Iwasawa triple $(A_M,b_MM,s_M)$ of~$G/M$ for which $s_M$ is smallest possible.
Note that if $N\le M$, then $s_N\ge s_M$, and assume for a contradiction that the set $S$ of all $s_M$ is unbounded.
Fix an open normal subgroup $N_1$ such that~$G/N_1$ is non-abelian and let $p^e$ be the exponent of $G/N_1$.
Let also $N_2$ be an open normal subgroup of $G$ such that $N_2\leq N_1$ and $s_{N_2}>e$. Using the natural epimorphism~\hbox{$G/N_2\rightarrow G/N_1$,} we see that $G/N_1$ should be abelian, a contradiction.

Therefore~$S$ has a largest element $s$.
Now, we can find an open normal subgroup~$L$ of $G$ containing $\overline{G'}$ and such that, for every $\overline{G'}\leq N\trianglelefteq_oG$ with $N\leq L$, the factor group~$G/N$ is the direct product of a cyclic group of order $\geq p^{s+1}$ and a finite group whose exponent is $\leq p^s$.
Let~$W/\overline{G'}$ be an infinite procyclic subgroup of $G/\overline{G'}$. Then~$G/W$ is a periodic group, because otherwise we would find in~$L$ an open normal subgroup~$N$ of $G$ such that~$G/N$ has two cyclic subgroups of order $p^{s+1}$ with trivial intersection.~Lem\-ma~4.3.7 of \cite{Zalesski} yields that $G/W$ has finite exponent, so also the subgroup~$T/\overline{G'}$ made by all elements of finite order of $G/\overline{G'}$ has finite exponent; in particular, $T/\overline{G'}$ is closed. Since $G/T$ is torsion-free, it is infinite and procyclic (see \cite{Zalesski}, Theorem 4.3.3).
This means we can write $G=X\ltimes T$ for some closed subgroup $X=\overline{\langle x\rangle}\simeq\mathbb{Z}_p$.

Fix~$N\trianglelefteq_o G$. Then $N\cap T$ is closed, and
$$
G/(N\cap T)\simeq X\ltimes\bigl(T/(N\cap T)\bigr).
$$
Since $T/N\cap T$ is finite, its centralizer $C$ in $X$ is a closed normal subgroup of $G$ with finite index in $X$.
Let $p^e=\operatorname{exp}(T/N\cap T)$ and, for any $f>e$, consider the factor group
$$
G/D(f)\simeq X/C^{p^f}\ltimes \bigl(T/(N\cap T)\bigr),
$$
where $D(f)=C^{p^f}(N\cap T)$.
By Theorem 2.3.1 of \cite{Schmidt}, the group $G/D(f)$ is a powerful $p$-group, and hence $\exp(G/D(f))=|X/C^{p^f}|$ (see Theorem 2.7 (iii) of \cite{DDMS99}).
In particular, the cyclic subgroups of maximal order of $G/D(f)$ are generated by elements of the form $x^naD(f)$, where $a\in T$ and~\hbox{$(n,p)=1$.} If any of them is normal, we have that $[X,T]\leq D(f)$, and $G/D(f)$ is abelian by~The\-o\-rem~2.3.15 of \cite{Schmidt}. If none of them is normal, we may use~The\-o\-rem~2.3.18 of \cite{Schmidt} to show that $T/(N\cap T)$ is abelian. Thus $T/(N\cap T)$ is abelian in any case. The arbitrariness of~$N$ yields that $T$ is abelian. Similar arguments prove that the action of $x$ on $T$ is precisely that described in point (iv).

Conversely, assume (3). Since the unrestricted direct product of modular lattices is modular, it is enough to show that groups of type (ii), (iii) and (iv) have a modular closed subgroup lattice.
Suppose first $G=X\ltimes A$, where $X$ and $A$ satisfy the properties described in (iv).
In this case, every finite quotient of $G$ is a modular $p$-group, so every pair of closed subgroups of $G$ is permutable by Lemma \ref{permute}; in particular, ${\mbox{\footnotesize$\vee$}}_c={\mbox{\footnotesize$\vee$}}$ and hence $\mathcal{L}_c(G)$ is modular, as required.

Assume now that $X$ and $A$ satisfy the properties in (iii) and let~\hbox{$H,K,L$} be closed subgroups of $G$ with $H\leq L$.
We need to show that
\begin{equation}
\label{equation}
\begin{array}{c}
\big(H\,{\mbox{\footnotesize$\vee$}}_c K\big)\,{\mbox{\footnotesize$\wedge$}}\,L=H\,{\mbox{\footnotesize$\vee$}}_c\, \big(K\,{\mbox{\footnotesize$\wedge$}}\, L\big).\tag{$\star$}
\end{array}
\end{equation}
If either $H$ or $K$ are contained in $A$, then (\ref{equation}) trivially holds.
Assume that neither $H$ nor~$K$ is contained in $A$; in particular, $L$ is not contained in $A$.
Write $$H=\overline{\langle x^{n_H}a_H\rangle}\ltimes A_H,\quad K=\overline{\langle x^{n_K}a_K\rangle}\ltimes A_K,\quad\textnormal{and}\quad L=\overline{\langle x^{n_L}a_L\rangle}\ltimes A_L,$$ where $A_H,A_K,A_L\leq A$, $a_H,a_K,a_L\in A$ and $n_H,n_K,n_L$ are non-zero $p$-adic exponents.
Let $$C=C_X(A)\cap\overline{\langle x^{n_H}a_H\rangle}\cap\overline{\langle x^{n_K}a_K\rangle}\cap\overline{\langle x^{n_L}a_L\rangle}.$$ Then $C$ has finite index in $X$ because each of the subgroups $\overline{\langle x^{n_H}a_H\rangle}, \overline{\langle x^{n_K}a_K\rangle}$, and $\overline{\langle x^{n_L}a_L\rangle}$ has finite index in $X\ltimes\langle a_H,a_K,a_L\rangle\in\mathcal{L}_c(G)$. Since $G/C$ is a $P^*$-group and $C$ is contained in $H\cap K\cap L$, we have $$\big(H\,{\mbox{\footnotesize$\vee$}}\,K\big)\,{\mbox{\footnotesize$\wedge$}}\,L=H\,{\mbox{\footnotesize$\vee$}}\,\big(K\,{\mbox{\footnotesize$\wedge$}}\,L\big).$$ We only have to show that $\langle H,K\rangle$ and $\langle H, K\cap L\rangle$ are closed. We prove more generally that the join of every two closed subgroups is closed. In fact, let $U$ and $V$ be two closed subgroups of $G$, and let $A_U=A\cap U$ and $A_V=A\cap V$. Then $A_U$ and $A_V$ are closed normal subgroups of $G$, so $A_UA_V$ is closed. As we did for $H, K$ and $L$, we can find a closed central subgroup $W$ of $G$ contained in $X\cap U\cap V$. Then $WA_UA_V$ is closed and has finite index in $\langle U,V\rangle$. This means that $\langle U,V\rangle$ is closed and we are done.

In order to complete the proof of the statement, we need to show that the join of any pair of closed subgroups of a locally finite and modular group is closed. This can be done by using the argument employed in the previous paragraph.
\end{proof}

\begin{cor}\label{cormodular}
Let $G$ be a profinite group. The following statements are equivalent:
\begin{itemize}
    \item[\textnormal{(1)}] Every pair of closed subgroups of $G$ is permutable;
    \item[\textnormal{(2)}] Every pair of open subgroups of $G$ is permutable;
    \item[\textnormal{(3)}] $G$ is an unrestricted direct product of coprime profinite groups each of which is of one of the following types: \begin{itemize}
        \item[\textnormal{(i)}] abelian;
        \item[\textnormal{(ii)}] locally finite, modular $p$-group for some prime $p$;
        \item[\textnormal{(iii)}] $X\ltimes A$, where $X\simeq\mathbb{Z}_p$ for some prime $p$, $A$ is an abelian pro-$p$ subgroup, and there are $k\in\mathbb{N}$ and $x\in X$  such that $X=\overline{\langle x\rangle}$ and $a^x=a^{1+p^k}$ for every $a\in A$; also $k\geq2$ if $p=2$.
    \end{itemize}
\end{itemize}
\end{cor}

Theorem \ref{modular} makes it also possible to obtain a description of a profinite group $G$ whose subgroup lattice is modular. In fact, $\mathcal{L}_o(G)$ is modular (being a sublattice of~$\mathcal{L}(G)$) and so $G$ satisfies part 3 of~The\-o\-rem~\ref{modular}. However, some further restrictions apply. In fact, it follows from~Lem\-ma~2.4.10 of~\cite{Schmidt} that profinite groups of types (iii) and (iv) cannot occur, and that $G$ must be the {\it restricted} direct product of groups of types (i) and (ii). Thus, the description of a profinite group with modular subgroup lattice coincides with that of an arbitrary locally finite group with modular subgroup lattice. Similar considerations can be made in case of a profinite group in which every pair of subgroups is permutable (here we use Corollary \ref{cormodular} instead of~Theo\-rem~\ref{modular}).

\medskip

We now turn to the study of modular elements of the closed subgroup lattice of a profinite group $G$ (note that every modular element of $\mathcal{L}_o(G)$ is easily seen to be also a modular element of $\mathcal{L}_c(G)$ by \cite{Schmidt}, Theorem 2.1.6).
The main result here is certainly Theorem \ref{5.1.14}, which deals with the structure of the factor group $G/M_G$, where $M$ is a modular element of $\mathcal{L}_c(G)$. Interesting consequences of this result are Corollary \ref{corcextmod}, which shows that the image of a modular element under an $o$-projectivity is always modular, and Theorem \ref{6.2.3}, a normality criterion for modular elements. Before proving~the main result, we prove a couple of related interesting facts concerning perfect profinite groups. The first one (Theorem \ref{5.2.6}) is an easy normality criterion, while the second one (Theorem \ref{5.4.9}) deals with the image of a closed normal subgroup under~\hbox{$o$-projectivities.}

\begin{theo}\label{5.2.6}
Let $G$ be a profinite group and let $M$ be a modular element of $\mathcal{L}_c(G)$ or $\mathcal{L}_o(G)$.
If $G$ or~$M$ is perfect, then $M$ is normal in $G$.
\end{theo}
\begin{proof}
Let $N$ be any open normal subgroup of $G$. It follows from Theorem 2.1.6 of \cite{Schmidt} that $MN/N$ is a modular element of $\mathcal{L}(G/N)$, so Corollary 5.2.6 of \cite{Schmidt} yields that $MN/N$ is normal in~$G/N$. The arbitrariness of $N$ shows that $M$ is normal in $G$.
\end{proof}

\medskip

Recall that if $\mathfrak{X}$ is a class of groups, the {\it $\mathfrak{X}$-residual} of a group $G$ is the intersection of all normal subgroups $N$ of $G$ such that $G/N$ belongs to $\mathfrak{X}$. In particular, if $\mathfrak{X}$ is the class of all soluble groups, we obtain the so-called {\it soluble residual} of the group $G$; and if $\mathfrak{X}$ is the class of all $p$-groups, we obtain the so-called {\it $p$-residual}.

\begin{theo}\label{5.4.9}
Let $G$ and $H$ be profinite groups and let $\varphi$ be an o-projectivity from $G$ onto $H$.
Let $N$ be a closed normal subgroup of $G$ such that $G/N$ is perfect. Then $N^\varphi\trianglelefteq H$.
\end{theo}
\begin{proof}
We need to show that $N^\varphi M/M$ is normal in $H/M$ for every $M\trianglelefteq_oH$. Let~$L$ be the normal core of $M_1=M^{\varphi^{-1}}$ in $G$. Since $G/N$ is perfect, we have $G=NR$, where~$R/L$ is the soluble residual of $G/L$.

We claim that $R^{\varphi}M/M$ is normal in $H/M$.
Observe that the lattice $\mathcal{L}(R/L)$ has no co-atom which is a modular element (see \cite{Schmidt}, Theorem 5.3.3, and Statement (1) at p.230), and since $M_1/L$ is a modular element of  $\mathcal{L}(G/L)$, it follows that also the lattice $$\big[(M_1\,{\mbox{\footnotesize$\vee$}}\,R)\,/\,M_1\big]\simeq\big[R/(R\,{\mbox{\footnotesize$\wedge$}}\, M_1)\big]$$ has no co-atom which is a modular element. This means that $R^\varphi M/M$ is perfect (see again \cite{Schmidt}, Theorem 5.3.3).
Moreover, $[G/R]$ satisfies the lattice condition described in Theorem 5.3.5 (d) of \cite{Schmidt}, and so also $\big[G/(M_1\,{\mbox{\footnotesize$\vee$}}\,R)\big]$ does.
Therefore $[H/R^\varphi M]$ satisfies the same condition, and $R^\varphi M/M$ is consequently the largest perfect subgroup of~$H/M$ (otherwise it would contradict Theorem 5.3.3 of \cite{Schmidt}).
This means that $R^\varphi M/M$ is the soluble residual of $H/M$; in particular,~$R^\varphi M/M$ is normal in $H/M$, as claimed.

Now, let $U=\big(N^\varphi M/M\big)^{H/M}$ and $V=\big(N^\varphi M/M\big)_{H/M}$. Since $N^\varphi M/M$ is a modular element of $\mathcal{L}(G/M)$, we have that $H/C_H(U/V)$ is supersoluble by Theorem 5.2.5 of \cite{Schmidt}, and hence that $R^\varphi$ is contained in $C_H(U/V)$. It follows that~$N^\varphi M/M$ is normalized by $R^\varphi$.
Since $H=N^\varphi R^\varphi$, we have that $N^\varphi M/M$ is normal in $H/M$.
\end{proof}

\begin{cor}
Let $G$ be a perfect profinite group. If $N$ is any closed normal subgroup of $G$, then $N^\varphi\trianglelefteq H$ for every $o$-projectivity $\varphi$ from $G$ onto a profinite group $H$.
\end{cor}

\medskip

Note that the arguments in the proof of Theorem \ref{5.4.9} can be easily employed to show that conditions such as prosolubility and prosupersolubility are invariant under~\hbox{$o$-pro}\-jectivities (see also \cite{Schmidt}, Theorem 5.3.7).

\medskip

Now, we prove our main result concerning modular elements of the closed subgroup lattice of a profinite group, but we first recall that a group~$G$ is a {\it $P$-group} if there is a prime $p$ and a cardinal number $n\geq2$ such that~$G$ is either elementary abelian of order $p^n$, or a semidirect product of an elementary abelian normal subgroup $A$ of order $p^{n-1}$ by a group of prime order $q\neq p$ which induces a non-trivial power automorphism on $A$; if $n$ is infinite, then $p^n$ and $p^{n-1}$ must be replaced by $n$. Recall also that a subgroup $H$ of an arbitrary group $G$ is {\it permutable} if $HK=KH$ for every~\hbox{$K\leq G$}; it is well known that permutable subgroups are ascendant, so in particular every permutable maximal subgroup must be normal.

\begin{theo}\label{5.1.14}
Let $G$ be a profinite group and let $M$ be a closed subgroup of $G$. The following properties are equivalent:
\begin{itemize}
    \item[\textnormal{(1)}] $M$ is a modular element of $\mathcal{L}_c(G)$;
    \item[\textnormal{(2)}] $M$ is a modular element of $\mathcal{L}_c\big(M\,{\mbox{\footnotesize$\vee$}}_c\,X\big)$ for every pro-$p$ procyclic subgroup $X$ of $G$;
    \item[\textnormal{(3)}] $MN/N$ is a modular element of $\mathcal{L}(G/N)$ for every $N\trianglelefteq_oG$;
    \item[\textnormal{(4)}] $G/M_G=\bigl(\operatorname{Cr}_{i\in I}S_i/M_G\bigr)\times T/M_G$, where $I$ is a \textnormal(possibly empty\textnormal) set of positive integers, the subgroups $S_i$ and $T$ are closed, and: \begin{itemize}
        \item[\textnormal{(a)}] $S_i/M_G$ is a non-abelian $P$-group for every $i\in I$;
        \item[\textnormal{(b)}] $\pi(S_i/M_G)\cap\pi^\star(T/M_G)=\emptyset$ for all $i\in I$ and $\pi(S_i/M_G)\cap\pi(S_j/M_G)=\emptyset$ for all $i\in I$ with $i\neq j$;
        \item[\textnormal{(c)}] $M/M_G=\bigl(\operatorname{Cr}_{i\in I} Q_i/M_G\bigr)\times (M\cap T)/M_G$, where $Q_i/M_G$ is a non-normal Sylow subgroup of $S_i/M_G$ for every $i\in I$, and $M\cap T$ permutes with every closed subgroup of $G$.
    \end{itemize}
\end{itemize}
\end{theo}
\begin{proof}
It is easy to see that (4) implies (1), and it is clear that (1) implies (2). 


Assume (2).  Let $N\trianglelefteq_oG$ and let $P/N$ be a cyclic subgroup of~$G/N$ of prime power order~$q^n$. Then there is a pro-$q$ procyclic subgroup $X$ of $G$ such that $P/N=XN/N$. By hypothesis,~$M$ is a modular element of $\mathcal{L}_c(U)$, where $U=M\,{\mbox{\footnotesize$\vee$}}_c\,X$, so $M(U\cap N)$ is a modular element of $\mathcal{L}_c(U/U\cap N)$; but $\mathcal{L}_c(U/U\cap N)$ is canonically isomorphic to~$\mathcal{L}_c(UN/N)$, and hence~$MN/N$ is modular in $\mathcal{L}(UN/N)=\mathcal{L}\big(\langle MN/N,P/N\rangle\big)$. Theorem~5.1.13 of~\cite{Schmidt} shows that $MN/N$ is a modular element of $G/N$, thus proving (3).

Assume (3). In order to prove (4), we may assume $M_G=\{1\}$. For any $N\trianglelefteq_oG$, let $L_N/N=(MN/N)_{G/N}$ --- Theorem 5.1.14 of \cite{Schmidt} shows that $G/L_N$ decomposes as in (4), although in this case, being $G/L_N$ finite, unrestricted direct products can be replaced by restricted ones. Let~$\pi$ be the set of all prime numbers~$p$ for which $G$ has an open normal subgroup $N_p$ such that $p$ divides the order of~$G/L_{N_p}$. Consider the partition $\sigma$ of $\pi$ constructed as follows: two distinct elements $p,q$ of $\pi$ are in the same element of $\sigma$ if one of the following conditions holds:
\begin{itemize}
    \item[(i)] there is an open normal subgroup $N$ of $G$ such that $$G/L_N=S_N/L_N\times T_N/L_N,$$ where 
    \begin{itemize}
        \item[(i1)] $S_N/L_N$ is a $P$-group whose order is divisible by $pq$,
        \item[(i2)] $\pi(S_N/L_N)\cap\pi(T_N/L_N)=\emptyset$, and
        \item[(i3)] $(MN/L_N)\cap(S_N/L_N)$ is a non-normal Sylow subgroup of $S_N/L_N$;
    \end{itemize}

    \item[(ii)] there is no open normal subgroup of $G$ satisfying condition (i) for the primes $p$ and $q$, and there is an open normal subgroup $N$ of $G$ such that $$G/L_N=S_N/L_N\times T_N/L_N,$$ where 
        \begin{itemize}
        \item[(i1)] $pq$ divides the order of $T_N/L_N$,
        \item[(i2)] $\pi(S_N/L_N)\cap\pi(T_N/L_N)=\emptyset,$ and
        \item[(i3)] $(MN/L_N)\cap(T_N/L_N)$ is a permutable subgroup of $G/L_N$.
        \end{itemize}

\end{itemize}

Note now that since $M_G\le L_N$ for every $N\trianglelefteq_oG$, we have $\bigcap_{N\trianglelefteq_o G}L_N=M_G=\{1\}$.
Hence, $G$ decomposes into the unrestricted direct product of coprime closed subgroups~$G_i$ such that $\pi(G_i)\in\sigma$. This means that we may assume $\pi(G)\in\sigma$, and that there exists an open normal subgroup $U$ of $G$ such that one of the following conditions holds for all open normal subgroups $N$ of $G$ which are contained in $U$:
\begin{itemize}
    \item[($\ast$)] $G/L_N$ is a non-abelian $P$-group such that $MN/L_N$ is a non-normal Sylow subgroup of $G/L_N$;
    \item[($\ast\ast$)] $MN/L_N$ is a permutable subgroup of $G/L_N$.
\end{itemize}

In case ($\ast\ast$), the subgroup $M$ permutes with all closed subgroups of $G$ (see Lem\-ma~\ref{permute}). Suppose ($\ast$) holds. Then $MN/L_N$ has prime order for any $N\trianglelefteq_oG$ with \hbox{$N\leq U$,} and we easily see that $L_KN/N=L_N/N$ for all $K\leq N\leq U$. It follows that $M$ has prime order, so every quotient $G/N$, with $N\trianglelefteq_oG$ and $N\leq U$, is a~\hbox{$P$-group}. Thus, $G$ has a closed normal subgroup $P$ which is an elementary abelian \hbox{$p$-group} for some prime~$p$, and $M$ acts on $P$ as a group of power automorphisms of prime order. This means that~$G$ is a $P$-group itself and the statement is proved.~\end{proof}

\begin{cor}
Let $G$ be a profinite group and let $M$ be a pro-$p$ subgroup of $G$. The following properties are equivalent:
\begin{itemize}
    \item[\textnormal{(1)}] $M$ is a modular element of $\mathcal{L}_c(G)$;
    \item[\textnormal{(2)}] either $M$ permutes with all closed subgroups of $G$, or $G/M_G=M^G/M_G\times K/M_G$, where $M^G/M_G$ is a non-abelian $P$-group which is coprime with $K/M_G$; moreover, both~$M^G$ and $K$ are closed subgroups of $G$.
\end{itemize}
\end{cor}

\begin{cor}
Let $G$ be a profinite group, and let $M$ be a maximal modular element of~$\mathcal{L}_c(G)$. Then either $M$ is normal in $G$ or $G/M_G$ is non-abelian of order $pq$ ($p$ and~$q$ primes).
\end{cor}

An interesting consequence of Theorem \ref{5.1.14} is concerned with the image of modular elements of the closed subgroup lattice under $o$-projectivities.

\begin{cor}\label{corcextmod}
Let $\varphi$ be the $c$-extension of an $o$-projectivity between the profinite groups $G$ and $H$, and let $M$ be a modular element of $\mathcal{L}_c(G)$. Then $M^\varphi$ is a modular element of $\mathcal{L}_c(H)$.
\end{cor}
\begin{proof}
In order to prove that $M^\varphi$ is a modular element of $\mathcal{L}_c(H)$ it suffices to show by The\-o\-rem~\ref{5.1.14} that $M^{\varphi}N/N$ is a modular element of $\mathcal{L}(G/N)$ for every $N\trianglelefteq_oG$.
Since $\varphi$ is the $c$-extension of an $o$-projectivity, we have
$$
\big(M\,{\mbox{\footnotesize$\vee$}}_c\,N^{\varphi^{-1}}\big)^\varphi=M^\varphi\,{\mbox{\footnotesize$\vee$}}_c\,N=M^\varphi N.
$$
Of course, $N^{\varphi^{-1}}$ is modular in $\mathcal{L}_o(G)$ and so it is modular in $\mathcal{L}_c(G)$.
Therefore $M\,{\mbox{\footnotesize$\vee$}}_c\,N^{\varphi^{-1}}$ is modular in the interval $\big[G/N^{\varphi^{-1}}\big]$ (see \cite{Schmidt}, Theorem 2.1.6), whence $M^\varphi N/N$ is a modular element of $\mathcal{L}(G/N)$. A final application of Theorem \ref{5.1.14} shows that $M^\varphi$ is modular in $G$.
\end{proof}

\medskip

As a consequence of Theorem \ref{5.1.14} we can also derive an interesting normality criterion for a modular element $M$ of the closed subgroup lattice of a profinite group $G$, but first we need some more information about the structure of $G/M_G$ when $[G/M]$ is a chain.

\begin{lem}\label{5.1.3}
Let $G$ be a profinite group. If $M$ is a modular element of $\mathcal{L}_c(G)$ such that $[G/M]$ is a chain, then either $G/M_G$ is a pro-$p$ group for some prime $p$, or $M$ is a maximal subgroup of $G$ and $G/M_G$ is non-abelian of order $pq$ for certain primes~$p$ and $q$. 
\end{lem}
\begin{proof}
Suppose first that for any $N\trianglelefteq_oG$, the factor group $G/L_N$ is a $p_N$-group, where $L_N/N=(MN/N)_{G/N}$ and $p_N$ is a prime.
Then, since $L_H\le L_N$ for any open normal subgroup $H$ of $G$ contained in $N$, and since $M_G=\bigcap_{N\trianglelefteq_o G}L_N$, it follows that $G/M_G$ is a pro-$p$ group for some prime $p$.


Assume therefore that there is some $N\trianglelefteq_oG$ such that the factor group $G/L_N$ is not a $p$-group. It follows from Theorem 2.1.6 of~\cite{Schmidt} that $MN/N$ is a modular element of~$\mathcal{L}_c(G/N)=\mathcal{L}(G/N)$, and of course $[G/MN]$ is still a chain. Thus, Lem\-ma~5.1.3 of \cite{Schmidt} yields that $G/L_N$ is non-abelian of order $pq$, where $p$ and $q$ are primes, and~$MN/N$ is maximal in $G/N$.
In particular $M$ is maximal in $G$, and, again, since $M_G=\bigcap_{N\trianglelefteq_o G}L_N$, it follows that $G/M_G$ is non-abelian or order $pq$, as desired.
\end{proof}

\medskip

A statement similar to that of Lemma \ref{5.1.3} holds if we replace $\mathcal{L}_c(G)$ by $\mathcal{L}_o(G)$, but in this case the result is an almost obvious consequence of \cite{Schmidt}, Lemma 5.1.3.

\begin{theo}\label{6.2.3}
Let $G$ be a profinite group, $M$ a modular element of $\mathcal{L}_c(G)$, and $g$ a non-trivial element of $G$ such that $\overline{\langle g\rangle}$ is torsion-free. If $M$ and $\overline{\langle g\rangle}$ are coprime, then $\overline{\langle g\rangle}$ normalizes~$M$.
\end{theo}
\begin{proof}
Of course, $\overline{\langle g\rangle}\cap M=\{1\}$, and we may assume $M_G=\{1\}$. Moreover, since $\overline{\langle g\rangle}$ is torsion free, we may also assume by Theorem \ref{5.1.14} that $M$ is permutable in $G$.

Let $H=\overline{\langle g\rangle}$ and write $X=MH$.
Since $M$ is modular in~$\mathcal{L}_c(G)$, we get $$[X/M]_c\simeq\mathcal{L}_c\big(\overline{\langle g\rangle}/\overline{\langle g\rangle}\cap M\big)=\mathcal{L}_c\big(\overline{\langle g\rangle}\big).$$ Thus, $M$ is the intersection of closed subgroups $L\in[X/M]_c$ such that $[X/L]_c$ is an infinite chain.
For any such subgroup $L$ there exists a prime $p\in\pi^\star(H)$ such that~\hbox{$L=MH_{p'}=:X_p$,} where $H_{p'}$ is the unique topological~\hbox{$p'$-Hall} subgroup of $H$. Since~$M$ is modular in $\mathcal{L}_c(X)$ and $[X/M]$ is distributive, it follows from~The\-o\-rem~2.1.6 of \cite{Schmidt} that $X_p$ is modular in $\mathcal{L}_c(X)$. Applying~Lem\-ma~\ref{5.1.3}, we have that $X/(X_p)_X$ is a pro-$p$ group. Since $p\not\in\pi^\star(X_p)$, it follows that $X_p=(X_p)_X$ is normal in $X$. Therefore~$M$ is the intersection of the normal subgroups $X_p$ and so it is normal in $X$.
\end{proof}

\medskip

We conclude this section with a couple of interesting (but easily provable) results on the topic. The first one is an easy application of Lemmas 5.1.10 and 5.1.11 of~\cite{Schmidt}, and our~Lem\-ma~\ref{permute} and Theorem \ref{5.1.14}.

\begin{theo}
Let $G$ be a profinite group and let $M$ be a pronilpotent subgroup of $G$. Then:
\begin{itemize}
    \item[\textnormal{(1)}] If $M$ permutes with every closed subgroup of $G$, then so does every topological Sylow subgroup of $M$.
    \item[\textnormal{(2)}] If $M$ is modular in $\mathcal{L}_c(G)$ and it is a topological Hall subgroup of $G$, then every Sylow subgroup of $M$ is modular in $G$.
\end{itemize}
\end{theo}

The second one can be proved by using~The\-o\-rem~5.1.7 of \cite{Schmidt} and similar arguments to those employed in the proof of Theorem~\ref{5.1.14} (one just needs to take into account the nilpotent residual of the modular element in every quotient by an open normal subgroup).

\begin{theo}
Let $G$ be a profinite group and let $M$ be a modular element of $\mathcal{L}_c(G)$. Then~$M/M_G$ is pronilpotent.
\end{theo}

\begin{flushleft}
\rule{8cm}{0.4pt}\\
\end{flushleft}

{
\sloppy
\noindent
Francesco de Giovanni

\noindent 
Dipartimento di Matematica e Applicazioni ``Renato Caccioppoli''

\noindent
Università degli Studi di Napoli Federico II

\noindent
Complesso Universitario Monte S. Angelo

\noindent
Via Cintia, Napoli (Italy)

\noindent
e-mail: degiovan@unina.it 

}

\bigskip
\bigskip

{
\sloppy
\noindent
Iker de las Heras

\noindent 
Zientzia eta Teknologia Fakultatea, Matematika Departamentua

\noindent
Euskal Herriko Unibertsitatea

\noindent
Sarriena Auzoa z/g, Leioa (Spain)

\noindent
e-mail: iker.delasheras@ehu.eus

}

\bigskip
\bigskip

{
\sloppy
\noindent
Marco Trombetti

\noindent 
Dipartimento di Matematica e Applicazioni ``Renato Caccioppoli''

\noindent
Università degli Studi di Napoli Federico II

\noindent
Complesso Universitario Monte S. Angelo

\noindent
Via Cintia, Napoli (Italy)

\noindent
e-mail: marco.trombetti@unina.it 

}


\begin{thebibliography}{99} 


\baselineskip 10pt
{

\bibitem{Acciarri}{\ssc C. Acciarri -- P. Shumyatsky}: ‘‘Profinite groups with restricted centralizers of $\pi$-elements’’, {\it Math. Z.} 301 (2022), 1039--1045.

\bibitem{liljan}{L. An}: ‘‘Groups whose Chermak-Delgado lattice is a subgroup lattice of an abelian group’’, {\it Canad. Math. Bull.} 66 (2023), 443--449.

\bibitem{Brandl}{\ssc R. Brandl}: ‘‘The Dilworth number of subgroup lattices’’, {\it Arch. Math. \textnormal(Basel\textnormal)} 50 (1988), 502--510.

\bibitem{DDMS99}{\ssc J.\,D.~Dixon -- M.\,P.\,F.~du~Sautoy -- A.~Mann --  D.~Segal}: ‘‘Analytic pro-p groups’’, Second edition, Cambridge Studies in Advanced Mathematics \textbf{61}, \textit{Cambridge University Press}, Cambridge (1999).


\bibitem{topft}{\ssc M. Ferrara -- M. Trombetti}: ‘‘The pro-norm of a profinite group’’, {\it Israel J. Math.} 254 (2023), 399--429.

\bibitem{FdGMT}{\ssc F. de Giovanni -- M. Trombetti}: ‘‘Large characteristic subgroups with modular subgroup lattice’’, {\it Arch. Math. \textnormal(Basel\textnormal)} 111 (2018), 123--128.

\bibitem{KhuShu}{\ssc E. Khukhro -- P. Shumyatsky}: ‘‘On profinite groups with automorphisms whose fixed points have countable Engel sinks’’, {\it Israel J. Math.} 247 (2022), 303--330.

\bibitem{Klopsch}{\ssc B. Klopsch -- M. Vannacci}: ‘‘Embedding properties of hereditarily just infinite profinite wreath products’’, {\it J. Algebra} 476 (2017), 297--310.

\bibitem{mann}{\ssc A. Mann}: ‘‘The exponents of central factor and commutator groups’’, {\it J. Group Theory} 10 (2007), 435--436.

\bibitem{Niksegal}{\ssc N. Nikolov -- D. Segal}: ``Finite index subgroups in profinite groups''. {\it C.R. Math. Acad. Sci. Paris} 337 (2003), 303--308.

\bibitem{ore}{\ssc O. Ore}: ‘‘Structures and group theory II’’, {\it Duke Math. J.} 4 (1938), 247--269.

\bibitem{Zalesski}{\ssc L. Ribes -- P. Zalesskii}: ‘‘Profinite Groups’’, {\it Springer}, Berlin (2000).

\bibitem{Rob72}{\ssc D.J.S. Robinson}: ‘‘Finiteness Conditions and Generalized Soluble Groups’’, {\it Springer}, Berlin (1972).

\bibitem{Rob82}{\ssc D.J.S. Robinson}: ‘‘A Course in the Theory of Groups’’, {\it Springer}, Berlin (1982).

\bibitem{ada}{\ssc A. Rottländer}: ‘‘Nachweis der Existenz nicht-isomorpher Gruppen von gleicher Situation der Untergruppen’’. {\it Math. Z.} 28 (1928), 641--653.

\bibitem{Schmidt}{\ssc R. Schmidt}: ‘‘Subgroup Lattices of Groups’’, {\it De Gruyter}, Berlin (1994).

\bibitem{Wilson}{\ssc J.S. Wilson}: ‘‘Profinite groups with few conjugacy classes of elements of infinite order’’, {\it Arch. Math. \textnormal(Basel\textnormal)} 120 (2023), 557–563.

\bibitem{Zacher}{\ssc G. Zacher}: ``Una caratterizzazione reticolare della finitezza dell'indice di un sottogruppo in un gruppo'', {\it Atti Accad. Naz. Lincei Rend. Cl. Sci. Fis. Mat. Nat.} 69 (1980), 317--323.

\bibitem{Zhang}{\ssc C. Zhang -- A.N. Skiba}: ‘‘On two sublattices of the subgroup lattice of a finite group’’, {\it J. Group Theory} 22 (2019), 1035--1047.



}
\end{thebibliography}
\end{document}